\documentclass[reqno]{amsart}

\usepackage[foot]{amsaddr}
\usepackage{graphicx,enumerate,nicefrac,bm}
\usepackage[dvipsnames]{xcolor}
\usepackage{cite}
\usepackage{amsmath, amssymb}
\usepackage{tikz}\usetikzlibrary{matrix}
\usepackage{a4wide}
\usepackage{subfig}
\usepackage{changes}
\usepackage{booktabs}
\usepackage{algorithm,algpseudocode,tabularx}

\makeatletter
\newcommand{\multiline}[1]{%
  \begin{tabularx}{\dimexpr\linewidth-\ALG@thistlm}[t]{@{}X@{}}
    #1
  \end{tabularx}
}
\makeatother

\newcommand{\E}{\mathsf{E}}

\newcommand{\x}{\bm{x}}

\newcommand{\dx}{\,\mathsf{d}\bm{x}}

\newcommand{\HS}{\mathrm{H}^1_0(\Omega)}

\newcommand{\Lom}[1]{\mathrm{L}^{#1}(\Omega)}
\renewcommand{\P}{\omega}
\newcommand{\Pref}{\widetilde{\P}}
\newcommand{\T}{\mathcal{T}}
\newcommand{\V}{\mathbb{V}}
\newcommand{\Vh}{\widehat{\mathbb{V}}}

\newcommand{\X}{\mathbb{X}}
\renewcommand{\H}{\mathbb{H}}

\newcommand{\poly}{\mathbb{P}}

\newcommand{\norm}[1]{\left\|#1\right\|}
\newcommand{\Lnorm}[1]{\left\|#1\right\|_{{\rm L}^2(\Omega)}}

\newcommand{\dprod}[1]{\langle #1\rangle}
\DeclareMathOperator{\Span}{span}

\newcommand{\tol}{\epsilon}
\newcommand{\incNk}{\mathrm{inc}_N^n}
\newcommand{\dENk}{\Delta\E_N^n}

\DeclareMathOperator*{\argmin}{arg\,min}

\newtheorem{theorem}{Theorem}[section]

\theoremstyle{definition}

\newtheorem{remark}[theorem]{Remark}

\title[Gradient Flow adaptive FEM for the GPE]{Gradient Flow Finite Element Discretizations with Energy-Based Adaptivity for the Gross-Pitaevskii Equation}

\author[P.~Heid]{Pascal Heid}
\email{pascal.heid@math.unibe.ch}

\author[B.~Stamm]{Benjamin Stamm}
\address{Center for Computational Engineering Science, Schinkelstr.~2, D-52062 Aachen, Germany}
\email{best@mathcces.rwth-aachen.de}

\author[T.~P.~Wihler]{Thomas P.~Wihler}
\address{Mathematics Institute, University of Bern, CH-3012 Bern, Switzerland}
\email{wihler@math.unibe.ch}

\thanks{The authors acknowledge the financial support of the Swiss National Science Foundation (SNF), Grant No. 200021\underline{\space}182524}

\begin{document}

\begin{abstract}
We present an effective adaptive procedure for the numerical approximation of the steady-state Gross--Pitaevskii equation. Our approach is solely based on energy minimization, and consists of a combination of gradient flow iterations and adaptive finite element mesh refinements. Numerical tests show that this strategy is able to provide highly accurate results, with optimal convergence rates with respect to the number of freedom.

\end{abstract}



\maketitle


\section{Introduction}

In quantum physics, \emph{Bose--Einstein condensates (BEC)} are important objects of study that feature various interesting properties including macroscopic quantum effects, superfluidity, and occurrence of quantum vortices (in the presence of a magnetic field). In order to model the steady states of BEC consisting of a collection of bosonic quantum particles, the time-independent \emph{Gross--Pitaevskii equation (GPE)} is widely used, see \cite{bose1924plancks,einstein1924quantentheorie,dalfovo1999theory}. It can be derived from the many-body Schr\"odinger equation with a given interaction potential in the limit of a large number of particles by applying a Hartree--Fock ansatz of a symmetric tensor product of a single-particle function (in contrast to a single determinant for fermions). Indeed, the Hartree--Fock ansatz becomes exact in the (dilute) mean-field limit (i.e.~under certain assumptions on the interaction between particles such as radial symmetry, repulsive and short-range), see \cite{lieb2001bosons,lewin2015mean} for some rigorous results. The GPE is a \emph{nonlinear eigenvalue problem} that represents the Euler-Lagrange equation of the Gross--Pitaevskii energy functional, given by
\begin{align} \label{eq:BEC}
 \E(v):=\int_\Omega\left( \frac{1}{2} |\nabla v|^2 + V(\x)|v|^2+\frac{\beta}{2} |v|^4\right) \dx,
\end{align}
under the following normalization constraint for the single particle functions~$v$:
\begin{equation}\label{eq:constraint}
\Lnorm{v}=1.
\end{equation}
Here, $\Omega \subset \mathbb{R}^d$, $d=\{1,2,3\}$, is a bounded, connected, and open set with Lipschitz boundary, $V\in L^\infty(\Omega)$ is a potential function with $V\geq 0$ almost everywhere, and $\beta\ge0$ is a constant. We note that the fourth-order term in~$\E$ (causing the associated eigenvalue to be nonlinear if~$\beta>0$) results from the interaction of particles. The global minimizer of~$\E$ under the constraint~\eqref{eq:constraint} is called the \emph{(normalized) ground state} of \eqref{eq:BEC}.   

\subsection*{State-of-the-art of numerical methods for BEC}

An overview of models and numerical schemes for the GPE is provided in \cite{bao2014mathematical}. In general terms, the minimizer of the Gross--Pitaevskii energy functional~\eqref{eq:BEC} can be found either 
\begin{enumerate}[(i)] 
\item by solving the corresponding Euler-Lagrange formulation, i.e.~the GPE,
\item or by direct minimization (under the normalization constraint~\eqref{eq:constraint}). 
\end{enumerate}

In the context of (i), classical spatial discretization approaches such as finite element methods~\cite{bao2003ground,gong2008finite,xie2016multigrid,raza2009energy}, finite difference schemes~\cite{chien2008two,raza2009energy}, Fourier methods~\cite{bao2003numerical,CancesChakirMaday:10}, or (pseudo-) spectral methods \cite{bao2005fourth,bao2006efficient,bao2009generalized} can be applied in a classical way. One possibility to deal with the nonlinearity occurring in the eigenvalue problem is to employ the Roothaan iteration scheme~\cite{cances2014perturbation}, also referred to as self-consistent field (SCF) iteration procedure. Alternatively, Newton's method~\cite{caliari2009minimisation} or adaptations to the inverse power method~\cite{jarlebring2014inverse} can be used. The idea of combining the iterative solution of the nonlinear eigenvalue problem and of mesh refinements in finite element discretizations has been addressed in 
\cite{cances2014perturbation,cances2017two,gong2008finite,chien2008two}.

Among the class of methods (ii) we point to the imaginary time method~\cite{bao2003numerical,antoine2014robust,antoine2014gpelab,bao2012mathematical} which, upon employing an imaginary time transformation, $t\to -\imath t$, is based on the observation that the time-dependent GPE relaxes to the ground state as time evolves. In addition, upon interpreting the imaginary time method as a steepest descent method for the energy functional~$\E$ from~\eqref{eq:BEC}, the closely related gradient flow approaches \cite{bao2010efficient,BaoDu:04,zeng2009efficiently,HenningPeterseim:18,bao2006efficient,danaila2010new,kazemi2010minimizing,raza2009energy} can be derived. Further numerical solution methods in the scope of~(ii) include a recently proposed preconditioned conjugate gradient method~\cite{antoine2017efficient}, direct energy minimization using symmetric properties~\cite{bao2003ground} (which simplify the constraint minimization problem), or the combination of gradient flows and Riemannian optimization~\cite{danaila2017computation}.

\subsection*{Contribution}

The aim of this work is to provide a numerical approximation procedure for the ground state of the Gross--Pitaevskii functional~\eqref{eq:BEC}, under the constraint~\eqref{eq:constraint}, which is based on a simultaneous interplay of gradient flow iterations and local adaptive finite element mesh refinements. This idea follows the recent developments on the (adaptive) \emph{iterative linearized Galerkin (ILG)} methodology~\cite{HeidWihler:18,HeidWihler:19,CongreveWihler:17,AmreinWihler:14,AmreinWihler:15,HoustonWihler:18}, whereby adaptive discretizations and iterative nonlinearity solvers are combined in an intertwined way; we also refer to the closely related works~\cite{ErnVohralik:13,El-AlaouiErnVohralik:11,BernardiDakroubMansourSayah:15,GarauMorinZuppa:11}.

A key building block of the numerical scheme to be presented in this paper concerns the decision of whether local mesh refinements or gradient flow iterations should be given preference on a given finite element discretization space. This is accomplished by monitoring the energy decay resulting from the gradient flow, and by performing a comparison to the energy loss caused by the latest mesh refinement. Then, depending on which effect is currently dominant, we either undertake another gradient flow step, or the mesh is refined adaptively. We emphasize that this is a very natural approach for the given problem since both the (conforming) finite element method and the recently proposed gradient flow method~\cite{HenningPeterseim:18} are both guaranteed to be energy-decreasing (owing to the variational principle). The proposed numerical method thereby generates a sequence of finite element approximations defined on adaptively refined spaces which provide a corresponding sequence of monotonically decreasing energies.

Proceeding along the energy minimization approach in~\cite{HoustonWihler:16}, the adaptive mesh-refinement strategy in this work is based on identifying a subset of elements in the mesh for which a local refinement will potentially provide a significant (local) contribution to the total energy decay. To this end, for each element in the mesh, we first apply a local gradient flow step on a locally refined patch; these computations, since local and independent, can be done in parallel and only involve very few degrees of freedom. Then, a D\"orfler marking strategy~\cite[Sec.~4.2]{Doerfler:96} selects the most promising elements for refinement.  The numerical tests illustrate optimal convergence rates in the number of unknowns for a variety of examples---both linear or nonlinear models, with smooth or irregular potentials, will be investigated.

\subsection*{Outline}
In Section~\ref{sec:GP} we present the framework of the Gross-Pitaevskii equation and motivate its associated gradient flow system. Moreover, Section~\ref{sec:FEM} is devoted to the finite element discretization and the adaptive mesh refinement strategy.
Furthermore, Section~\ref{sec:numerics} presents various numerical tests in 2D. Finally, we add some concluding remarks in Section~\ref{sec:conclusions}.

\section{The Gross--Pitaevskii equation and gradient flows}\label{sec:GP}

\subsection{Nonlinear eigenvalue formulation}

We observe that the energy functional $\E$ from~\eqref{eq:BEC} is Fr\'{e}chet differentiable on the Sobolev space~$\H:=\HS$ (defined to be the space of all functions with weak gradient in~$\Lom{2}^d$ and zero trace along the boundary~$\partial\Omega$); indeed, a simple calculation reveals that 
\begin{align} \label{eq:BECderivative}
 \dprod{\E'(v),w}= \int_\Omega \left(\nabla v \cdot \nabla w+ 2V(\x)vw+2 \beta |v|^2vw\right)\dx,\qquad v,w\in\H,
\end{align}
where $\dprod{\cdot,\cdot}$ signifies the dual product in $\H^\star \times \H$. The Euler-Lagrange formulation of the constrained minimization problem 
\[
 \argmin_{v \in S_\H} \E(v),
\]
with $S_\H:=\{v \in \H: \|v\|_{{\rm L}^2(\Omega)}=1\}$ signifying the $\Lom{2}$-unit sphere in~$\H$, is given by 
\begin{align} \label{eq:GP}
 v \in \H: \qquad \dprod{\E'(v),w}=\lambda (v,w)_{{\rm L}^2(\Omega)} \qquad \forall w \in \H.
\end{align}
Here, the scalar $\lambda$ takes the role of a Lagrange multiplier corresponding to the norm constraint~\eqref{eq:constraint}. The \emph{nonlinear} eigenvalue problem~\eqref{eq:GP} is called the \emph{Gross--Pitaevskii equation (GPE)}. If $v \in S_\H$ is an eigenfunction of \eqref{eq:GP} associated with an eigenvalue $\lambda$, then we note that
\begin{align} \label{eq:EVenergy}
 \lambda=\dprod{\E'(v),v}=2 \E(v)+\beta \|v\|^4_{{\rm L}^4(\Omega)}>0.
\end{align}

Given that~$V\ge 0$ (almost everywhere in~$\Omega$) and~$\beta\ge 0$, the Gross--Pitaevskii eigenvalue problem~\eqref{eq:GP} has a unique ($\Lom{2}$-normalized) positive eigenfunction~$u_{\mathrm{GS}}>0$ which is the ground state of the Bose--Einstein condensate~\eqref{eq:BEC}, see~\cite[Lem.~5.4]{HenningPeterseim:18}; in particular, $u_{\mathrm{GS}}$ is an eigenfunction to the minimal (and simple) eigenvalue, signified by~$\lambda_{\mathrm{GS}}$ of~\eqref{eq:GP}, see \cite{CancesChakirMaday:10}. 

\subsection{Continuous gradient flow}

The ground state $u_{\mathrm{GS}}$ will be determined iteratively. To this end, we employ the \emph{projected gradient flow} approach proposed in \cite{HenningPeterseim:18}. One of the key ideas is to introduce a weighted energy inner product on~$\H\times\H$: For fixed $z \in \H$, we let
\begin{align} \label{eq:innerproduct}
 a_z(v,w):=\int_\Omega \left(\nabla v \cdot \nabla w+2 V(\bm{x})vw+2\beta |z|^2vw\right) \dx,\qquad v,w\in\H.
\end{align}
Owing to the Riesz representation theorem, for any~$v\in\H$, it exists a unique $\mathsf{G}_z(v)\in\H$ such that
\begin{align} \label{eq:green}
 a_z(\mathsf{G}_z(v),w)=(v,w)_{\Lom{2}} \qquad \forall w \in \H.
\end{align}
If $z\neq 0$, we notice that $a_z(\mathsf{G}_z(z),z)=\|z\|^2_{\Lom{2}}>0$, and therefore~$\mathsf{G}_z(z)\neq 0$. Hence, we may consider the linear mapping
\begin{equation}\label{eq:P}
\mathsf{P}_{z}:\H \to \mathbb{T}_{z},\qquad
\mathsf{P}_{z}(v)=v-\frac{(z,v)_{{\rm L}^2(\Omega)}}{a_z(\mathsf{G}_z(z),\mathsf{G}_z(z))}\mathsf{G}_z(z),
\end{equation}
where
\[
\mathbb{T}_{z}:=\{w \in \H: (z,w)_{\Lom{2}}=0\}; 
\]
for~$z=0$, we have~$\mathbb{T}_{0}=\H$, and $\mathsf{P}_{0}$ is the identity map. Using~\eqref{eq:green}, it is fairly elementary to verify that
\begin{align} \label{eq:orthproj}
a_z(v-\mathsf{P}_{z}(v),w)=0\qquad\forall v\in\H,\,\forall w\in \mathbb{T}_{z},
\end{align}
i.e.~$\mathsf{P}_{z}$ is the orthogonal projection onto the (tangential plane) $\mathbb{T}_{z}$ with respect to the $a_z$-inner product.

Based on the above definitions, we are now ready to present the gradient flow induced by the inner product from~\eqref{eq:innerproduct}. More precisely, we consider a trajectory~$u:\,[0,\infty)\to\H$ which, for a given initial value~$u(0)=u_0\in S_\H$, follows the dynamical system
\begin{align} \label{eq:continuousflow}
\dot{u}(t)=-\mathsf{P}_{u(t)}(u(t)), \qquad t > 0.
\end{align}
The existence of a solution~$u$ has been discussed in~\cite[Sec.~3.1 \& 3.2]{HenningPeterseim:18}. 

In order to motivate the gradient flow system~\eqref{eq:continuousflow}, we briefly revisit some arguments from~\cite[Proof of Theorem~3.2]{HenningPeterseim:18}:
\begin{enumerate}[(i)]
\setlength{\itemsep}{5pt}
\item For any $v,z\in\H$, we define~$\nabla_z\E(v)\in\H$ to be the Riesz representative of the steepest descend direction~$\E'(v)$ at~$v$ with respect to the $a_z$-inner product, i.e.
\[
a_z(\nabla_z\E(v),w)=\dprod{\E'(v),w}\qquad\forall w\in\H.
\]
Then, using~\eqref{eq:continuousflow}, employing~\eqref{eq:orthproj} respectively~\eqref{eq:innerproduct}, and recalling~\eqref{eq:BECderivative}, for any~$t>0$, we notice that
\begin{align*}
a_{u(t)}(\dot u(t),w)
&=a_{u(t)}(-\mathsf{P}_{u(t)}(u(t)),w)\\
&=a_{u(t)}(-u(t),\mathsf{P}_{u(t)}(w))\\
&=\dprod{-\E'(u(t)),\mathsf{P}_{u(t)}(w)}\\
&=a_{u(t)}(-\nabla_{u(t)}\E(u(t)),\mathsf{P}_{u(t)}(w))\\
&=a_{u(t)}(\mathsf{P}_{u(t)}(-\nabla_{u(t)}\E(u(t))),w),
\end{align*}
for any~$w\in\H$. This implies that the gradient flow from~\eqref{eq:continuousflow} follows the orthogonal projection of the steepest descend direction $-\nabla_{u(t)}\E(u(t))$, for~$t>0$, onto the tangential plane~$\mathbb{T}_{u(t)}$.
\item Observe that
\begin{equation}\label{eq:udot}
\dot u(t)\in \mathbb{T}_{u(t)}\qquad \forall t>0.
\end{equation} 
Hence,
\begin{align*}
 \frac{\mathsf{d}}{\mathsf{dt}} \norm{u(t)}_{\Lom{2}}^2=2 (u(t),\dot u(t))_{\Lom{2}}=0\qquad\forall t>0.
\end{align*}
In particular, it follows that $\norm{u(t)}_{\Lom{2}}=\norm{u(0)}_{\Lom{2}}=1$. This means that the gradient flow stays on the sphere~$S_\H$ for any~$t>0$, i.e., physically speaking, it is mass preserving; cf.~\cite[Lem.~3.3]{HenningPeterseim:18}. 
\item Using again~\eqref{eq:udot} and invoking~\eqref{eq:orthproj}, for any~$t>0$, it holds that
\begin{equation}\label{eq:Edecay}
\frac{\mathsf{d}}{\mathsf{dt}}\E(u(t))
=\dprod{\E'(u(t)),\dot u(t)}
=a_{u(t)}(u(t),\dot u(t))
=-a_{u(t)}(-\mathsf{P}_{u(t)}(u(t)),\dot u(t)).
\end{equation}
Thus, implementing~\eqref{eq:continuousflow} yields
\[
\dprod{\E'(u(t)),\dot u(t)}=-a_{u(t)}(\dot u(t),\dot u(t))\le 0\qquad\forall t>0,
\]
i.e. the energy~$\E(t)$ is monotone decreasing as~$t\to\infty$; cf.~\cite[Lem.~3.3]{HenningPeterseim:18}. 
\item Since~$\E(v)$ is nonnegative for any~$v\in\H$, the monotonicity property from~(iii) implies that there is~$\E^\star\ge 0$ with~$\lim_{t\to\infty}\E(u(t))=\E^\star<\infty$. Hence, applying~\eqref{eq:Edecay}, we obtain the identity
\[
0\le\int_0^\infty a_{u(\tau)}(\dot u(\tau),\dot u(\tau))\,\mathsf{d}\tau=\E(u_0)-\E^\star<\infty,
\]
which implies that $\int_0^\infty\|\nabla \dot u(\tau)\|^2_{\Lom{2}}\,\mathsf{d}\tau$ is bounded. It follows that~$u(t)$ has a limit $u^\star\in \H$ as~$t\to\infty$ (with~$\Lnorm{u^\star}=1$, cf.~(ii)), and~$\mathsf{P}_{u^\star}(u^\star)=0$ upon taking the limit in~\eqref{eq:continuousflow}. Consequently, we deduce from~\eqref{eq:P} that $u^\star=\lambda^\star\mathsf{G}_{u^\star}(u^\star)$, where
\[
\lambda^\star=\frac{\|u^\star\|^2_{\Lom{2}}}{a_{u^\star}(\mathsf{G}_{u^\star}(u^\star),\mathsf{G}_{u^\star}(u^\star))}.
\]
Therefore, exploiting~\eqref{eq:green}, for all~$w\in\H$, we infer that
\begin{equation}
	\label{eq:NonLinEVPlimit}
\dprod{\E'(u^\star),w}
=a_{u^\star}(u^\star,w)
=\lambda^\star a_{u^\star}(\mathsf{G}_{u^\star}(u^\star),w)
=\lambda^\star(u^\star,w)_{\Lom{2}},
\end{equation}
i.e. $u^\star\not\equiv0$ is an eigenfunction for the GPE~\eqref{eq:GP} to the eigenvalue~$\lambda^\star$.
We note that $u^\star$ can be any eigenfunction satisfying~\eqref{eq:NonLinEVPlimit}, respectively \eqref{eq:GP}; in particular, it is not necessarily the ground state.
\end{enumerate}


\subsection{Discrete gradient flow}

For the purpose of computing an approximation of the continuous gradient flow trajectory from~\eqref{eq:continuousflow}, we use the forward Euler discretization method. For a given initial value $u^0 \in S_\H$ this yields a sequence of functions~$\{u^n\}_{n\ge 0}\subset S_\H$ which, for~$n\ge0$, is defined by
\begin{subequations}\label{eq:GFiteration}
\begin{align} 
  u^{n+1}&=\frac{\widehat{u}^{n+1}}{\Lnorm{\widehat{u}^{n+1}}},\label{eq:GFiteration1}
  \intertext{where}
  \widehat{u}^{n+1}&
  =u^n-\tau_n\mathsf{P}_{u^n}(u^n)
  =(1-\tau_n)u^n+\frac{\tau_n}{a_{u^n}(\mathsf{G}_{u^n}(u^n),\mathsf{G}_{u^n}(u^n))} \mathsf{G}_{u^n}(u^n).\label{eq:GFiteration2}
\end{align}
\end{subequations}
Here, $\{\tau_n\}_{n\ge 0}$ is a sequence of positive (discrete) time steps that is assumed uniformly bounded from above and below with bounds~$\tau_{\max}$ and~$\tau_{\min}$, respectively, such that $0 < \tau_{\min} \leq \tau_n \leq \tau_{\max} < \infty$. The scheme~\eqref{eq:GFiteration} is called \emph{discrete gradient flow iteration (GFI)}.

\begin{remark}\label{rem:tau}
Provided that the maximal time step $\tau_{\max}$ is sufficiently small, it can be seen that the GFI scheme~\eqref{eq:GFiteration} yields guranteed energy reduction. More precisely, there exists $0<\tau_{\max}=\mathcal{O}(\min\{\beta^{-1},\E(u^0)^{-\nicefrac{1}{2}}\})$ such that for all $\tau_n \leq \tau_{\max}<2$ it holds that $\E(u^{n+1}) \leq \E(\widehat{u}^{n+1}) \leq \E(u^n)$, see~\cite[Lemma~4.7]{HenningPeterseim:18}. Moreover, if $\tau_n \leq 1$ for all $n \geq 0$ then, for any starting value $u^0 \in S_\H$ with $u^0 \geq 0$, the (full) sequence $\{u^n\}_{n \geq 0}$ generated by the GFI~\eqref{eq:GFiteration} satisfies~$u^n\ge 0$ for all~$n\ge 0$ and converges strongly in $\H$ to the unique positive ground state $u_{\mathrm{GS}}$, see~\cite[Theorem 5.1]{HenningPeterseim:18}.
\end{remark}

\section{Adaptive gradient flow finite element discretization}\label{sec:FEM}

We now focus on the adaptive spatial discretization of the gradient flow iteration scheme~\eqref{eq:GFiteration}. 

\subsection{Finite element discretization}

Consider a sequence of conforming and shape-regular partitions $\{\T_N\}_{N\in\mathbb{N}}$ of the domain~$\Omega$ into simplicial elements~$\T_N=\{\kappa\}_{\kappa\in\T_N}$ (i.e.~triangles for $d=2$ and tetrahedra for~$d=3$). Moreover, for a (fixed) polynomial degree~$p\in\mathbb{N}$ and any subset~$\omega \subset\T_N$, we introduce the finite element space
\[
\V(\omega)=\left\{v\in\H:\,v|_\kappa\in\poly_p(\kappa), \kappa\in\omega,\,v|_{\Omega\setminus\omega}=0\right\},
\]
with~$\poly_p(\kappa)$ signifying the (local) space of all polynomials of maximal total degree~$p$ on~$\kappa$, $\kappa\in\T_N$. Furthermore, similarly as before, we denote by
\[
S_{\V(\omega)}=\left\{v\in\V(\omega):\,\|v\|_{\Lom{2}}=1\right\}
\]
the $\Lom{2}$-unit sphere in~$\V(\omega)$. In the sequel, we apply the notations $\X_N:=\V(\mathcal{T}_N)$ and $S_N:=S_{\V(\mathcal{T}_N)}$.
We further denote by $\E_{N}:=\E|_{\X_N}$ the restriction of the energy functional $\E$ from \eqref{eq:BEC} to the Galerkin space $\X_N$. Then, due to the compactness of $S_N$, it exists a minimizer $u_{N} \in \X_N$ of $\E_N$, i.e. $\E(u_N)=\min_{v \in S_{N}} \E(v)$, with $(u_{N},1)_{{\rm L}^2(\Omega)} \geq 0$.
It is not known, however, if $u_N$ is unique (up to the sign), see, e.g.,~\cite{CancesChakirMaday:10}.  
Furthermore, if $\{\X_N\}_{N\in\mathbb{N}}$ is a dense family of finite element subspaces of $\H$, then any sequence of minimizers $u_{N} \in S_{N}$ of $\E_{N}$ with $(u_{N},1)_{{\rm L}^2(\Omega)} \geq 0$  converges in $\mathrm{H}^1(\Omega)$ to the ground state $u_{\mathrm{GS}}$ of $\E$; we refer to \cite[Theorem 1]{CancesChakirMaday:10} or \cite[Theorem 3.1]{Zhou:04}.

%

\subsection{Discrete GFI}

Let us define the (space) discrete version of the gradient flow iteration~\eqref{eq:GFiteration} on a  finite element subspace $\X_N \subset \H$. For $u \in \X_N$ we denote by $\mathsf{G}^N (u) \in \X_N$ the unique solution of 
\begin{align}\label{eq:RieszN}
 	a_{u}(\mathsf{G}^N(u),v)=(u,v)_{{\rm L}^2(\Omega)} \qquad \forall v \in \X_N;
\end{align}
cf.~\eqref{eq:green}. 
For given $u\in\X_N$, note that the computation of $\mathsf{G}^N(u)$ is a standard linear source problem; it can be solved by any linear solver at the disposal of the user. 
Then, for~$n\ge 0$, the space discrete GFI is given by
\begin{subequations}\label{eq:discreteGF}
\begin{align} 
 	u_N^{n+1}
 	&=\frac{\widehat{u}_N^{n+1}}{\norm{\widehat{u}_N^{n+1}}_{\Lom{2}}}, \label{eq:discreteGF1}
 	\intertext{where}
 	\widehat{u}_N^{n+1}&=(1-\tau_{N}^n) u_N^n+\frac{\tau_{N}^n}{a_{u_N^n}(\mathsf{G}^N(u_N^n),\mathsf{G}^N(u_N^n))}\mathsf{G}^N(u_N^n),\label{eq:discreteGF2}
\end{align}
\end{subequations}
with a sequence of discrete time steps~$\{\tau_N^n\}_{n \geq 0}$ as in \eqref{eq:GFiteration}. 

\begin{remark}
Consider a fixed mesh $\T_N$ and associated approximation space~$\X_N$. Let $\{u_N^n\}_{n \geq 0}\subset\X_N$ be the sequence generated by the discrete GFI~\eqref{eq:discreteGF} with some initial value $u_N^0 \in S_{N}$. If $\tau_N^n \leq \tau_{\max}$, with $\tau_{\max}$ as in Remark~\ref{rem:tau}, for all $n \geq 0$, then the corresponding energies are strictly monotone decreasing, and there exists a limit energy $\E_N^\star=\lim_{n \to \infty} \E(u_N^n)$. Furthermore, up to subsequences, we have $u_N^n \to u_N^\star$ strongly in $\H$, where $u_N^\star \in S_N$, with $\E(u_N^\star)=\E_N^\star$, is a discrete eigenfunction of the corresponding GPE, i.e.~there is $\lambda_N^\star$ so that
\begin{equation}\label{eq:evpd}
	a_{u_N^\star}(u_N^\star,v)=\lambda_N^\star(u_N^\star,v)_{{\rm L}^2(\Omega)} \qquad \forall v \in \X_N.
\end{equation}
We refer to \cite[Corollary 4.11]{HenningPeterseim:18} for details.
\end{remark}

\begin{remark}\label{rem:tauN}
In practical computations, in order to guarantee a positive energy decay in each iteration step, we propose the time step strategy within~\eqref{eq:discreteGF} given by
\begin{align*} 
 \tau_N^{n} = \max\left\{2^{-m}:\, \E(u_N^{n+1}(2^{-m})) < \E(u_N^n),\,m\ge 0\right\},\qquad  n \geq 0.
\end{align*}
where, for~$0<s\le 1$, we write $u_N^{n+1}(s)$ to denote the output of the discrete GFI \eqref{eq:discreteGF} based on the time step $\tau_N^n=s$ and on the previous approximation~$u_N^n$. We observed in several examples that for the choice $\tau_N^0=1$, i.e.~using $m=0$ above, no time correction was needed; we also refer to~\cite[Remark 4.8]{HenningPeterseim:18} for a discussion of the fixed time step $\tau=1$. For that reason, and for the sake of keeping the computational cost minimal, we fix the time step $\tau=1$ in the local GFI from Algorithm~\ref{alg:ref} below. We still use, however, the time step strategy for the global GFI in Algorithm~\ref{alg:adaptive}.
\end{remark}

\subsection{Local energy decay and adaptive mesh refinements}
\label{ssec:LocEn}

\begin{figure}
\begin{center}
\begin{tabular}{cc}
\includegraphics[scale=0.2]{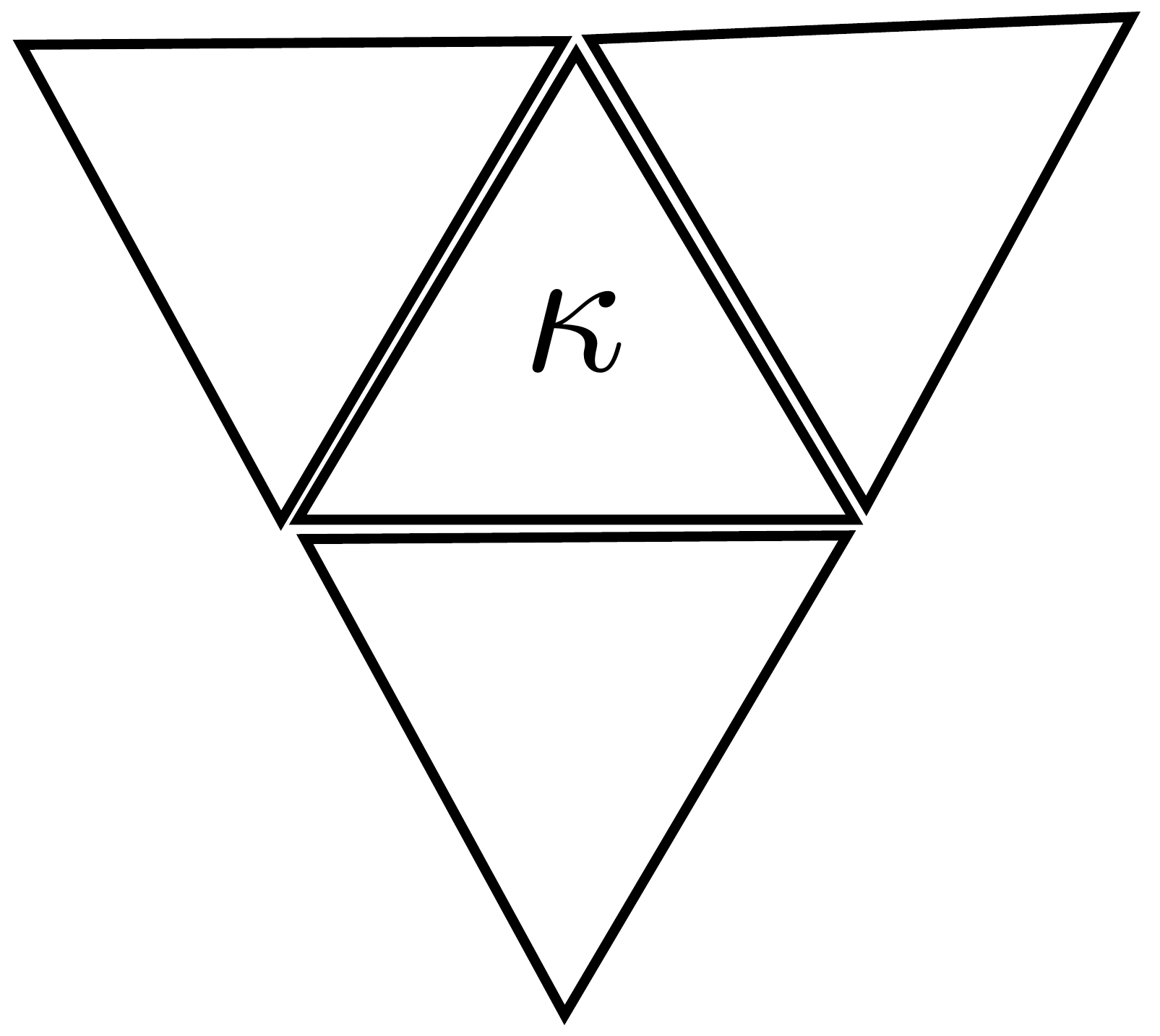} &
\includegraphics[scale=0.2]{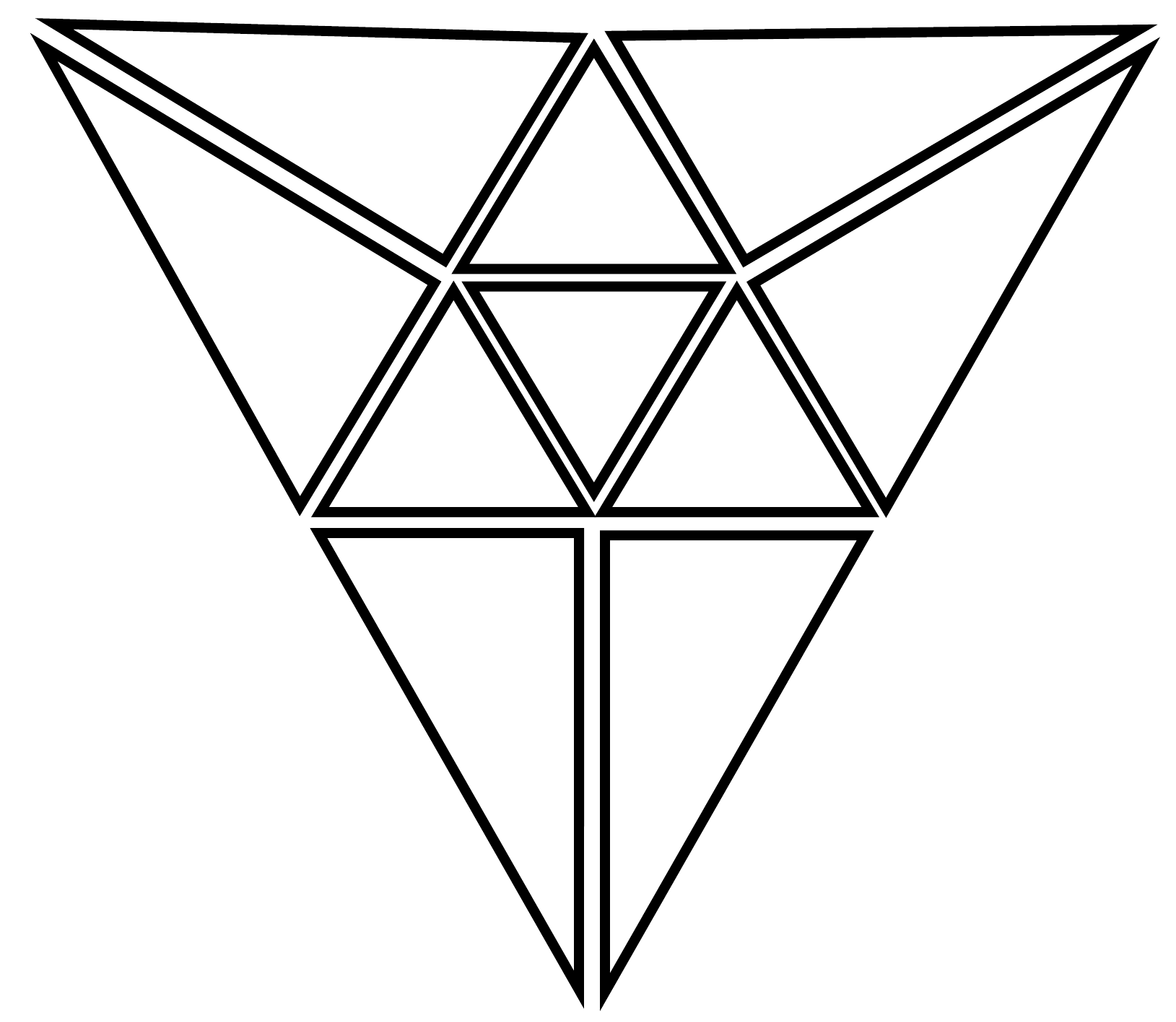}
\end{tabular}
\end{center}
\caption{Local element patches associated to a triangular element~$\kappa$. Left: Mesh patch~$\P_\kappa$ consisting of the element $\kappa$ and its face-neighbours. Right: Modified patch $\Pref_\kappa$ constructed based on red-refining $\kappa$ and on green-refining its neighbours.}
\label{fig:ref}
\end{figure}

For any element~$\kappa\in\T_N$ we consider the open patch~$\P_\kappa$ comprising of $\kappa$ and its immediate face-wise neighbours. Moreover, given $\kappa \in \T_N$, we define the modified patch~$\Pref_\kappa$ by uniformly (red) refining the element $\kappa$ into a (fixed) number of subelements; here, we assume that the introduction of any hanging nodes in $\P_\kappa$ is removed by doing (e.g. green) refinements, see Figure~\ref{fig:ref}. 
We remark that the notions of \emph{red and green refinements} refer to standard element subdivision techniques in automatic mesh adaptation; further details can be found, e.g., in~\cite[\textsection{}4.10.2.2]{LarsonBengzon:13} (or \cite{BankShermanWeiser:83}, where \emph{red refinement} is termed \emph{bisection-type mesh refinement}).

We consider basis functions~$\{\xi^1_\kappa,\ldots,\xi^{m_\kappa}_\kappa\}$ of the locally supported space~$\V(\Pref_\kappa)$. Furthermore, for any given~$v\in\V(\T_N)$, we introduce the extended space
\[
\Vh(\Pref_\kappa;v):=\Span\{\xi^1_\kappa,\ldots,\xi^{m_{\kappa}}_\kappa,v\}.
\] 
Suppose we have found an accurate approximation $u_N^n \in \X_N$ of the discrete GPE~\eqref{eq:evpd}, for some~$n\ge 0$. Then, by performing one local discrete GFI-step in $\Vh(\Pref_\kappa;u_N^n) \subset \H$ we obtain a new local approximation, denoted by $\widetilde u_{N,\kappa}^n \in \Vh(\Pref_\kappa;u_N^n)$, with $\norm{\widetilde u_{N,\kappa}^n}_{\Lom{2}}=1$. We emphasize that $\Vh(\Pref_\kappa;u_N^n)$ has a small dimension, and hence the discrete GFI~\eqref{eq:discreteGF} based on $\Vh(\Pref_\kappa;u_N^n)$ entails hardly any computational cost (for instance, for dimension $d=2$ and polynomial degree~$p=1$, the dimension of the locally refined space $\Vh(\Pref_\kappa;u_N^n)$, cf.~Figure~\ref{fig:ref}, is typically 3 or 4). 

By modus operandi of the discrete GFI~\eqref{eq:discreteGF}, the above construction leads in general, see Remark~\ref{rem:tauN}, to the (local) energy decay
\begin{align} \label{eq:localenergydecay}
 -\Delta \E_N^n(\kappa):=\E(\widetilde u_{N,\kappa}^n)-\E(u_N^n) \leq  0,
\end{align}
for all $\kappa \in \mathcal{T}_N$. The value $\Delta \E_N^n(\kappa)$ indicates the potential energy reduction due to a refinement of the element $\kappa$. This observation motivates the energy-based adaptive mesh refinement procedure outlined in Algorithm~\ref{alg:ref}.   

\begin{algorithm}
\caption{Energy-based adaptive mesh refinement}
\label{alg:ref}
\begin{algorithmic}[1]
\State Prescribe a mesh refinement parameter~$\theta\in(0,1)$. 
\State Input a finite element mesh~$\T_N$, and an $\Lom{2}$-normalized finite element function~$u^{n}_N \in \X_N$ for some~$n\ge 1$.
\For {all elements $\kappa\in\T_N$}
	\State\multiline{Perform one discrete GFI-step in the low-dimensional space $\Vh(\Pref_\kappa;u_N^n)$ to obtain a potentially improved local approximation~$\widetilde u_{N,\kappa}^n$.}  
	\State Compute the local energy decay~$\Delta \E_N^n(\kappa)$ from~\eqref{eq:localenergydecay}.    
\EndFor
\State \textsc{Mark} a subset ~$\mathcal{K} \subset \mathcal{T}_N$ of minimal cardinality which fulfils the D\"orfler marking criterion
\[
\sum_{\kappa \in \mathcal{K}} \Delta \E_N^n(\kappa) \geq \theta \sum_{\kappa \in \mathcal{T}_N} \Delta \E_N^n(\kappa).
\]
\State \textsc{Refine} all elements in~$\mathcal{K}$ for the sake of generating a new mesh~$\T_{N+1}$.
\end{algorithmic}
\end{algorithm}

\subsection{Adaptive strategy}

From a practical viewpoint, once the discrete GFI approximation from \eqref{eq:discreteGF} is close to a solution of~\eqref{eq:evpd}, on a given finite element space, we expect that any further GFI steps will no longer contribute an essential decay to the energy in~\eqref{eq:GFIincrement}. In this case, in order to further reduce the energy, we need to enrich the finite element space appropriately. More specifically, for~$N \ge 1$, suppose that we have performed a reasonable number~$n \ge 1$ (possibly depending on~$N$) of GFI-iterations~\eqref{eq:discreteGF} in $\X_{N-1}$. Consider now a (hierarchically) refined mesh~$\T_{N}$ of~$\T_{N-1}$, for example, obtained by Algorithm~\ref{alg:ref}. Then we may embed the final guess $u^{n}_{N-1}\in \X_{N-1}$ on the previous space into the enriched finite element space~$\X_{N}$ in order to obtain an initial guess on the refined mesh~$\T_{N}$:
\[
	u^0_{N}:=u^{n}_{N-1}\in \X_{N}.
\]
For each GFI-iteration $n$ we monitor two quantities. Firstly, we introduce the increment on each iteration given by
\begin{align*} 
	\incNk := \E(u_N^{n-1}) -  \E(u_N^n),\qquad n\ge 1.
\end{align*}
Secondly, we compare $\incNk$ to the energy loss as compared to the previous mesh refinement, i.e.
\begin{align} \label{eq:GFIincrement}
	\dENk := \E(u_N^{0}) -  \E(u_N^n),\qquad n\ge1.
\end{align}
We stop the iteration for~$n\ge 1$ as soon as $\incNk$ becomes small compared to $\dENk$, i.e.~once there is no notable benefit (relatively speaking) in performing any more discrete GFI steps on the current space $\X_{N}$. Specifically, for~$n\ge 1$, this is expressed by the bound
\[
	\incNk \le \gamma\dENk,
\]
for some parameter $0<\gamma<1$. 
%
%
We implement this procedure in Algorithm~\ref{alg:adaptive}.

\begin{algorithm}
\caption{Adaptive finite element gradient flow procedure}
\label{alg:adaptive}
\begin{algorithmic}[1]
\State Prescribe the three parameters~$\theta,\gamma\in(0,1)$, and $0<\epsilon\ll 1$.
\State Choose a sufficiently fine initial mesh~$\T_0$, and an initial guess~$u^0_0\in S_\H$ with $u^0_0 \geq 0$. Set~$N:=0$.
\Loop 
	\State Set $n:=1$, perform one discrete GFI-step in $\X_N$ to obtain $u_{N}^{1}$.
	\State Compute the indicator $\mathrm{inc}_N^1$ (which equals~$\Delta\E_N^1$). 
	\While {$\incNk > \gamma\dENk$} 
		\State Update $n\gets n+1$.
		\State Perform one GFI-step in $\X_N$ to obtain $u_{N}^{n}$ (starting from $u_N^{n-1}$).
		\State Compute the indicators $\incNk$ and $\dENk$.
	\EndWhile
	\If {$\dENk > \tol \E(u_N^n)$}
		\State \multiline{\textsc{Mark} and adaptively \textsc{refine} the mesh~$\T_N$ using Algorithm~\ref{alg:ref} to generate a new mesh~$\T_{N+1}$.}
		\State Define $u_{N+1}^0:=u_N^n\in \X_{N+1}$ by canonical embedding~$\X_N\hookrightarrow\X_{N+1}$.
		\State Update~$N\gets N+1$.
	\Else 
	\State \Return $u^n_N$.
	\EndIf
\EndLoop
\end{algorithmic}
\end{algorithm}

\subsection{Computational complexity}
We comment on the computational cost of one loop occurring in Algorithm~\ref{alg:adaptive}. This is essentially comprised of a number of GFI steps on a given Galerkin space $\X_N$, and on one adaptive (local) mesh refinement (using Algorithm~\ref{alg:ref}).
\begin{itemize}
\item The cost of one GF iteration is dominated by the computation of the Riesz-representative $\mathsf{G}^N(u_N^n)$ from~\eqref{eq:RieszN}, and, thereby, depends primarily on the dimension $\textrm{dim}(\mathbb X_N)$ of the finite element space $\mathbb X_N$. Specifically, the solution of~\eqref{eq:RieszN} (for given $u^n_N$) amounts to a computational work that scales like $\mathcal O(\textrm{dim}(\mathbb X_N)^\alpha)$, with a parameter $\alpha\ge 1$ depending on the linear solver employed.

\item The local finite element space $\widehat{\mathbb V}(\widetilde\omega_\kappa;u_N^n)$ contains $m_\kappa$ local and one global basis function, namely~$u_N^n$, which is the same for all elements $\kappa$. The mesh-refinement procedure contains one GF iteration as described in~\eqref{eq:discreteGF}, however, on the local space $\widehat{\mathbb V}(\widetilde\omega_\kappa;u_N^n)$. This, in turn, requires the solution of a linear system involving $m_\kappa$ local and one global degrees of freedom. In particular, only one entry of the corresponding matrix requires a global integration, 
which can be computed element-by-element (and, therefore, in parallel); all other matrix entries are represented by local integrals. Finally, note that the entry requiring global integration is \emph{the same for all local GF iterations and all elements $\kappa$}, and, thus, needs to be computed only once within each step of the loop in Algorithm~\ref{alg:adaptive}. 

A similar observation holds for the computation of the local energy decays. In fact, for $u_N^n \in \X_N$, denoting, as in Algorithm~\ref{alg:ref}, the locally improved approximation by $\widetilde{u}_{N,\kappa}^n$, we have the linear combination
\[
\widetilde{u}_{N,\kappa}^n=\mu u_N^n + \sum_{i=1}^{m_\kappa} \mu_i \xi_\kappa^i,
\] 
for suitable $\mu,\mu_i \in \mathbb{R}$, $i \in \{1,\dotsc,m_{\kappa}\}$. Then, it is elementary to verify that 
\[
\E(\widetilde{u}_N^n)=\E(\mu u_N^n)-\E_{\widetilde{w}_\kappa}(\mu u_N^n)+\E_{\widetilde{w}_\kappa}(\widetilde{u}_{N,\kappa}^n),
\] 
where $\E_{\widetilde\omega_\kappa}(\cdot)$ represents the contribution from the elements comprising the patch $\widetilde\omega_\kappa$ to the total energy $\E(\cdot)$. We emphasize that the second and third terms on the right-hand side only require the computation of a local integral. Moreover, for the term $\E(\mu u_N^n)$ we observe that
$
\E(\mu u_N^n)=:\mu^2 \E^n_{N,1} + \mu^4 \E^n_{N,2},
$
with
\begin{align*}
\E^n_{N,1}
&= \int_\Omega\left(\frac{1}{2} |\nabla u_N^n|^2 + V(\x)|u_N^n|^2\right) \dx,&
\E^n_{N,2}&=\int_\Omega \frac{\beta}{2} |u_N^n|^4 \dx.
\end{align*}
Thus, for each patch, the energy $\E(\mu u_N^n)$ is a combination of the above two global integrals, which need to be computed once only, and can be split into elementwise contributions.

Altogether, the computational cost for the local GF iterations and corresponding local energy decays therefore scales linearly with the number of elements in the mesh, i.e. with $\textrm{dim}(\mathbb X_N)$, and can be performed fully in parallel. The same observation applies to the element marking and local mesh refinement procedure, which compares to the evaluation of standard residual based a posteriori error estimators.
\end{itemize}
In summary, we find that the computational work for one loop in Algorithm~\ref{alg:adaptive} scales with $\mathcal O(\textrm{dim}(\mathbb X_N)^\alpha)$, where $\alpha\ge 1$ depends on the global finite element solver. In particular, our scheme exhibits a similar complexity as standard adaptive finite element discretization procedures for linear problems (provided that the number of GFI steps remain reasonably modest). 

\section{Numerical Experiments}\label{sec:numerics}

We apply Algorithm~\ref{alg:adaptive} for some numerical computations in two space dimensions, i.e.~$d=2$, with Cartesian coordinates denoted by $\bm{x}=(x,y) \in \mathbb{R}^2$. In all examples, we choose the initial guess $u_0^0 \in S_0$ such that $u_0^0(\bm{x})=c$ for any node $\bm{x}$ in the interior of the corresponding (coarse and uniform) initial mesh $\mathcal{T}_0$, where $c>0$ is the appropriate constant to fulfil the norm constraint~\eqref{eq:constraint}; we remark that numerical experiments based on other positive initial guesses (not presented here) have resulted in analogous convergence rates of the minimal energy approximation, thereby indicating a certain robustness with respect to the starting values. Moreover, we set $\theta=0.5$, as well as $\gamma=0.1$ and $\tol=10^{-8}$ in Algorithms~\ref{alg:ref} and~\ref{alg:adaptive}, respectively. Even if the stopping criterion in Line~12 of Algorithm~\ref{alg:adaptive} may not be satisfied, for the purpose of our tests, we stop the computations once the number of degrees of freedom (i.e.~the dimension of the finite element space~$\X_N$) exceeds~$10^6$.

\subsection{Laplace EVP on $L$-shaped domain} \label{sec:Lshaped}
\begin{figure}
	\hfill
	\includegraphics[width=0.49\textwidth]{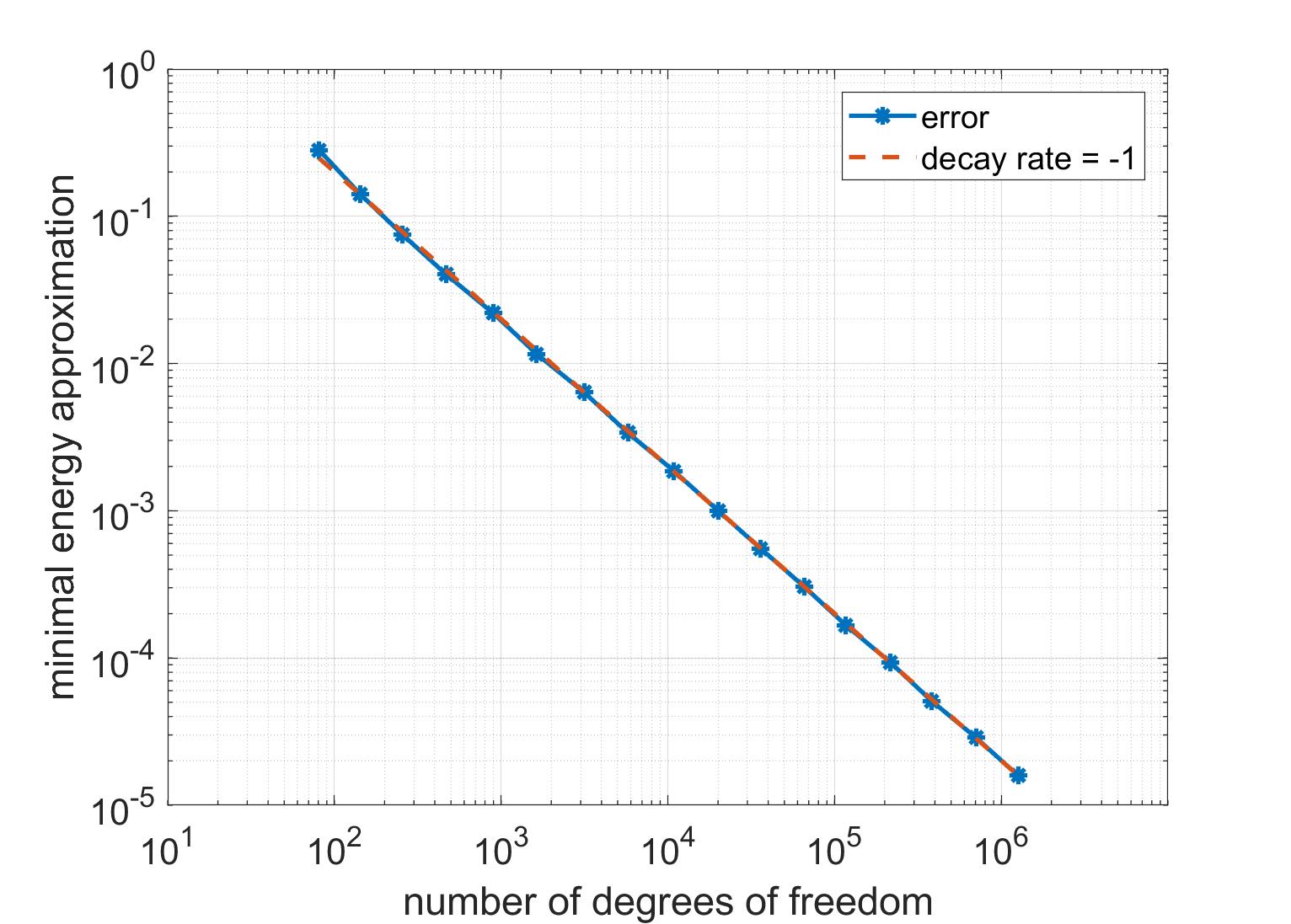}
	\hfill
	\includegraphics[width=0.49\textwidth]{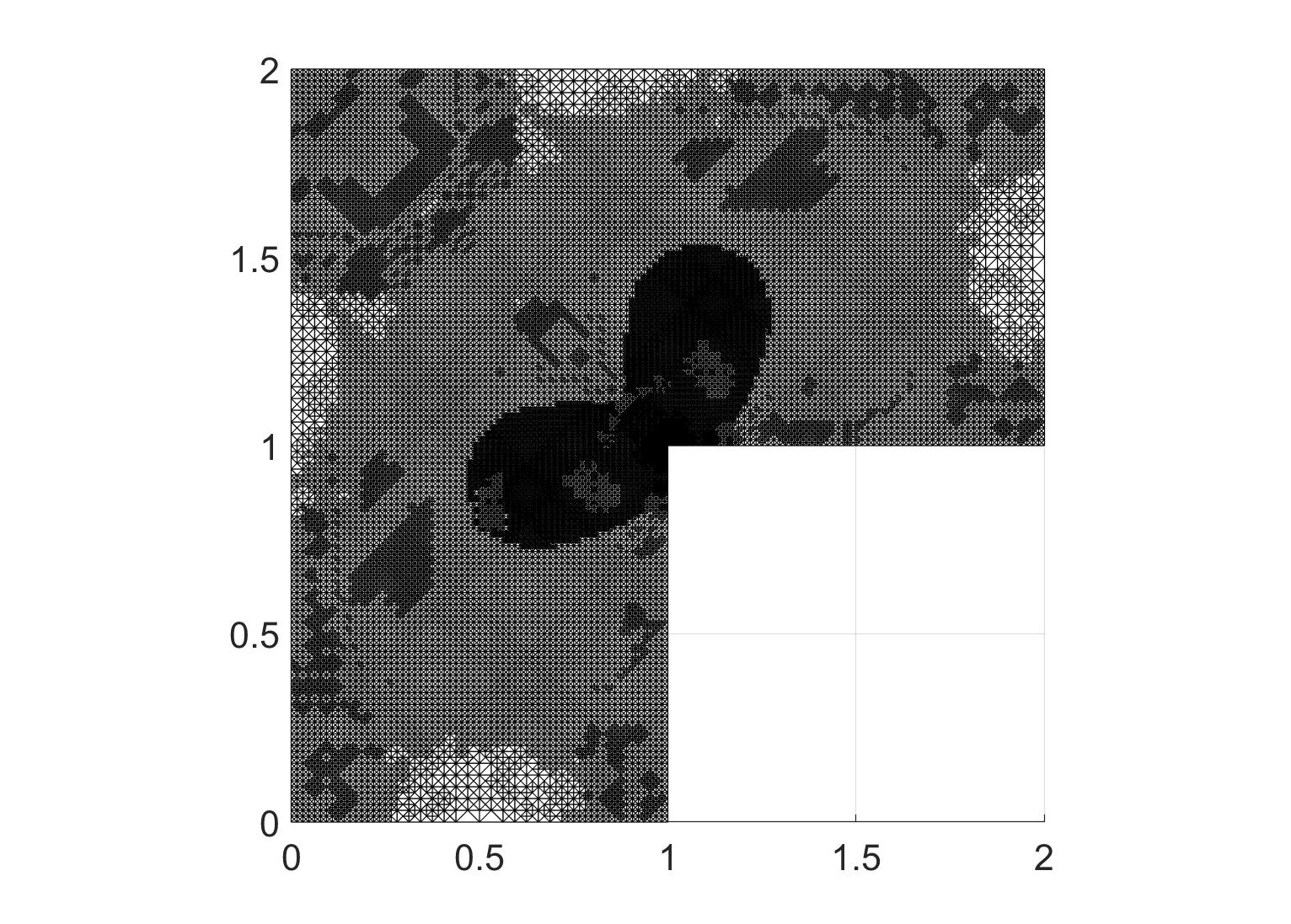}
	\hfill
	\caption{Experiment~\ref{sec:Lshaped} Left: Convergence plot for the ground state energy. Right: Adaptively refined mesh after 10 refinements.}
	\label{fig:Lshape}
\end{figure}

We begin by testing our algorithm for the Laplace eigenvalue problem, which is to find $u\in\H=\HS$ and $\lambda>0$ such that
\begin{align} \label{eq:laplace}
-\Delta u= \lambda u\qquad\text{in }\Omega;
 \end{align}
here, $\Omega=(0,2)^2 \setminus [1,2] \times [0,1]$ is an $L$-shaped domain. This problem is of interest since the eigenfunction to the lowest eigenvalue $\lambda_1$, i.e.~the ground state, has a singularity at the re-entrant corner point $(1,1)$. From~\cite{FoxMoler:67} it is known that $9.6397238 \leq \lambda_1 \leq 9.6397239$. Taking the lower bound (divided by 2, cf.~\eqref{eq:EVenergy} for $\beta=0$) as reference value, Figure~\ref{fig:Lshape} demonstrates the optimal convergence rate for the minimal energy approximation.
 
\subsection{Linear GPE with smooth potential $V$} \label{sec:y9}

\begin{figure}
	\hfill
	\includegraphics[width=0.49\textwidth]{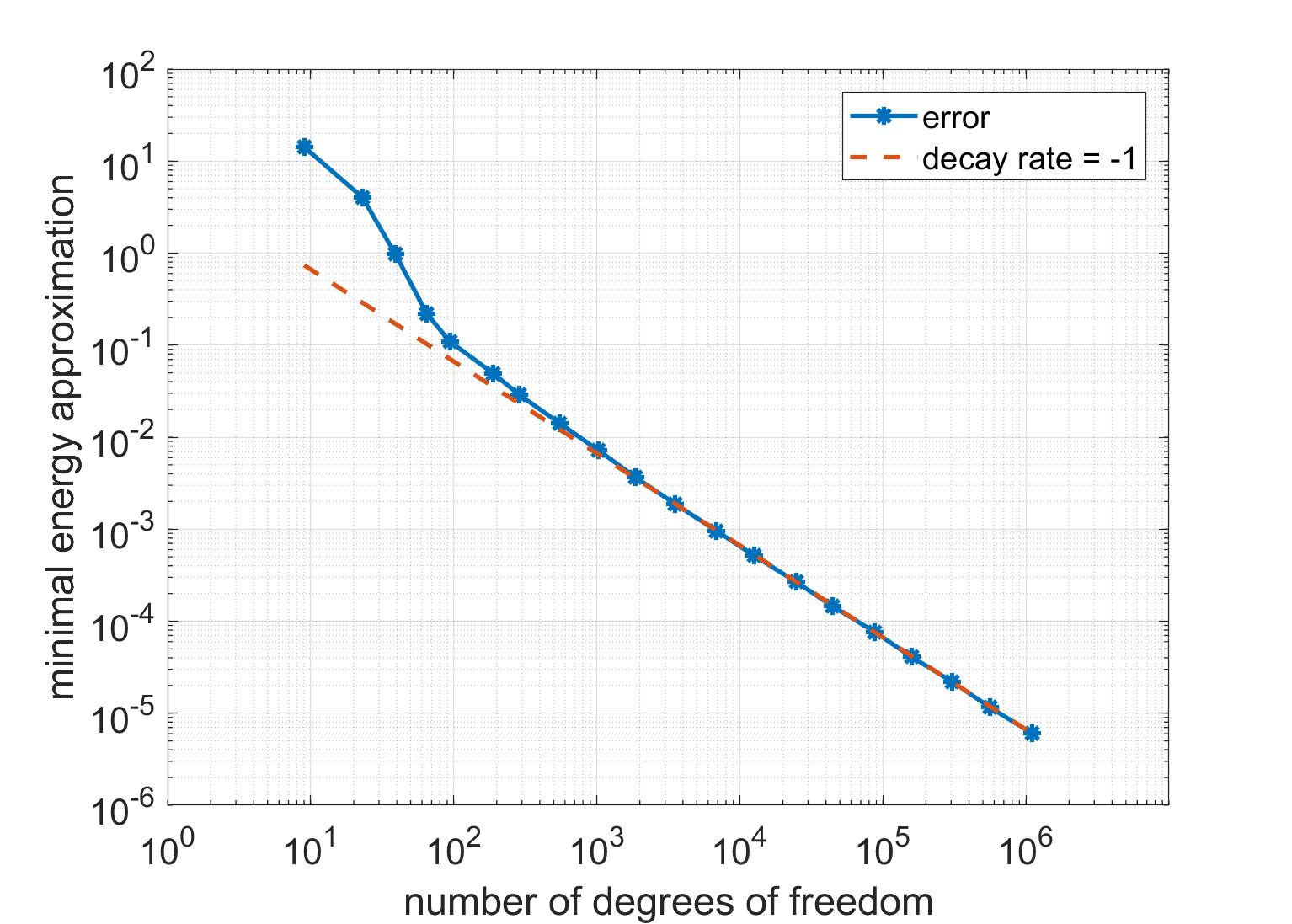}
	\hfill
  	\includegraphics[width=0.49\textwidth]{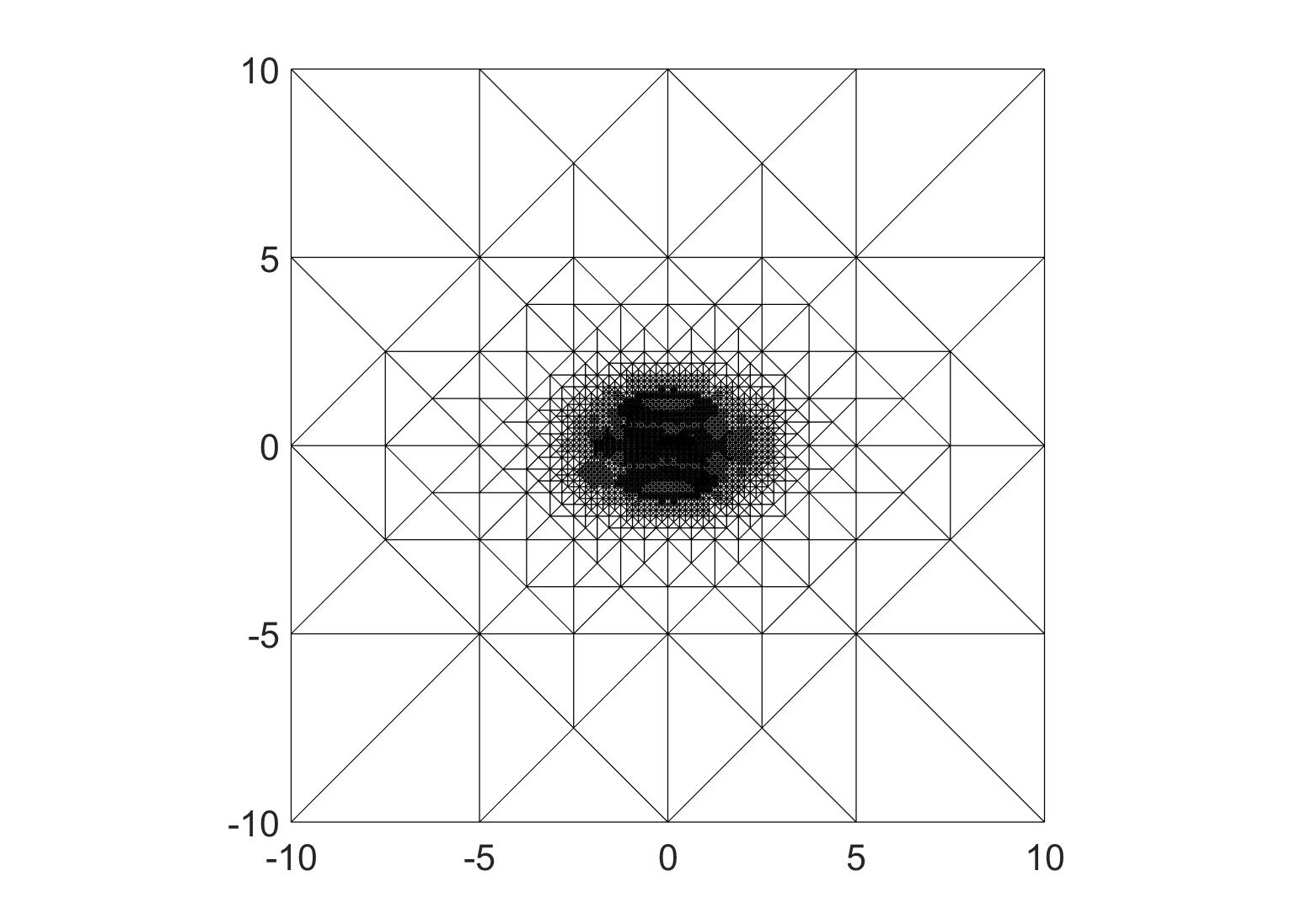}
	\hfill
	\caption{Experiment~\ref{sec:y9}. Left: Convergence plot for the energy of the ground state. Right: Adaptively refined mesh after 10 refinements.}
\label{fig:y9}
\end{figure}

We consider the case where~$\beta=0$ in~\eqref{eq:BEC}, and~$V(x,y)=\nicefrac{1}{2}(x^2+9y^2)$ is a smooth function:
\begin{align*}
 \E(u)=\frac{1}{2} \int_\Omega \left(|\nabla u|^2+\left(x^2+9y^2\right)|u|^2\right) \dx.
\end{align*}
Note that the associated eigenvalue problem~\eqref{eq:GP} is linear for this example. It is known that, for $\Omega=\mathbb{R}^2$, the energy of the ground state is given by $\E(u_{\mathrm{GS}})=2$, see, e.g.,~\cite{bao2014mathematical}. We note that the mass of $u_{\mathrm{GS}}$ is essentially concentrated in a vicinity of the origin $\bm{0}$ due to the global minimum of~$V$ at~$\bm{0}$. Therefore, the restriction to the bounded domain $\Omega:=(-10,10)^2 \subset \mathbb{R}^2$, which we use in our computations, has almost no effect on the minimal value of the ground state energy. Figure~\ref{fig:y9} (left) illustrates that the approximation of the energy converges with optimal rate in terms of the numbers of degrees of freedom. In addition, we see that the mesh is mainly refined around the origin, where the mass of the ground state $u_{\mathrm{GS}}$ is concentrated, see Figure~\ref{fig:y9} (right). These results underline that the proposed adaptive gradient flow procedure effectively detects the local behaviour of the model.

\subsection{Linear GPE with singular potential $V$.} \label{sec:alpha1}

\begin{figure}
	\hfill
	\includegraphics[width=0.49\textwidth]{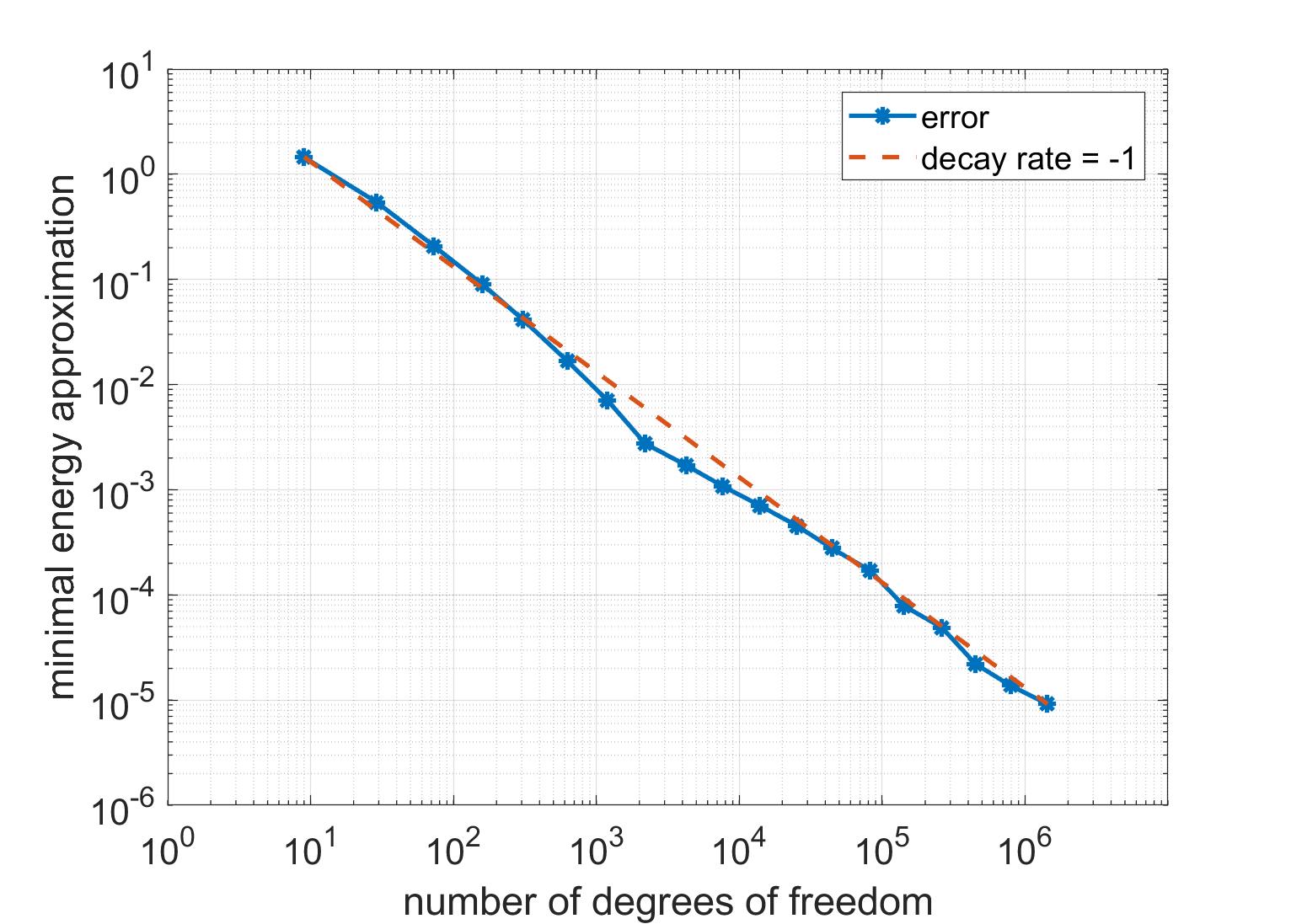}
	\hfill
	\includegraphics[width=0.49\textwidth]{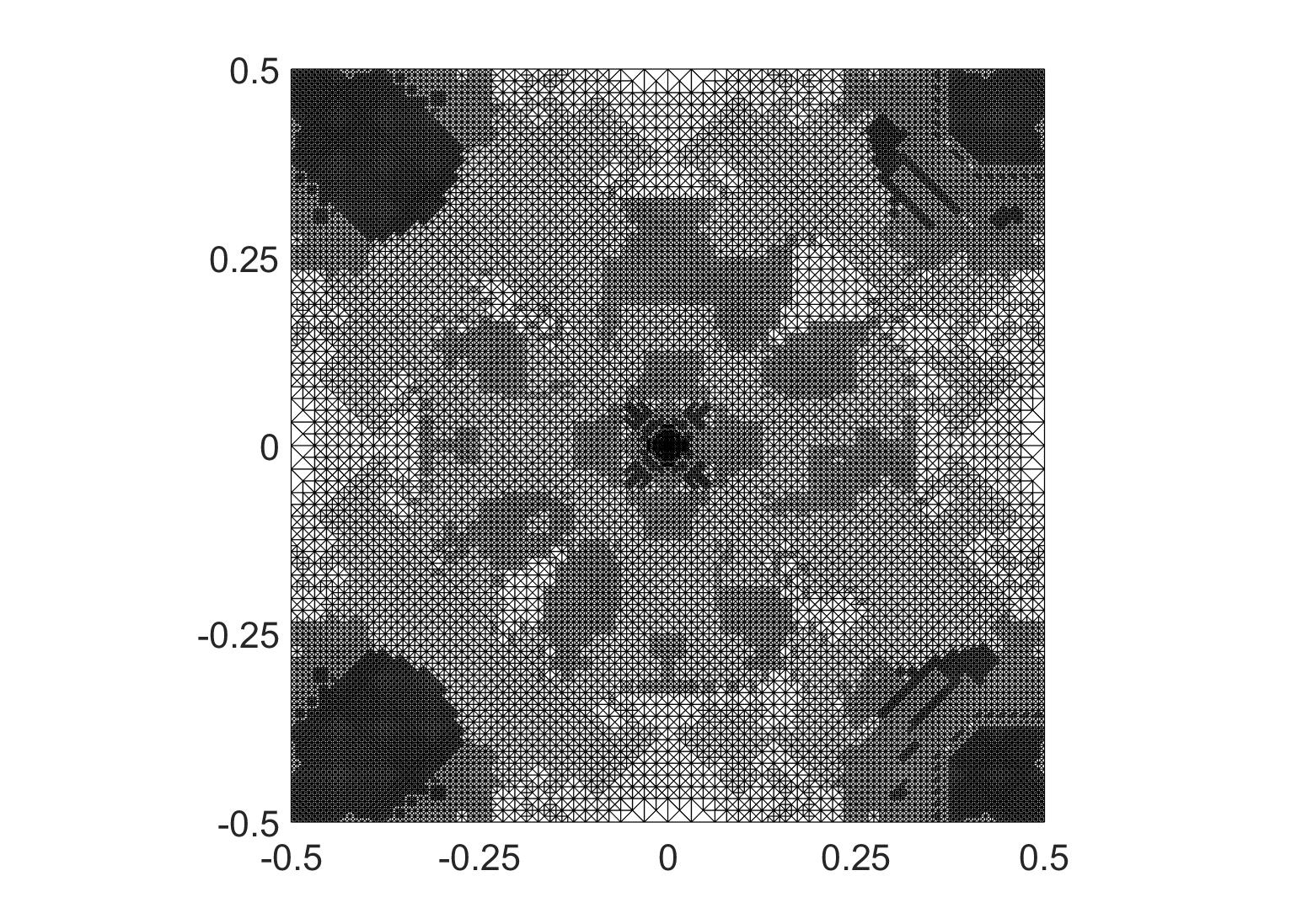}
	\hfill
	\caption{Experiment~\ref{sec:alpha1}. Left: Convergence plot for the energy of the ground state. Right: Adaptively refined mesh after 11 refinements.}
	\label{fig:alpha1}
\end{figure}

We perform another experiment with~$\beta=0$ in~\eqref{eq:BEC}, i.e.~the associated eigenvalue problem~\eqref{eq:GP} is again linear. In contrast to the previous example, however, we consider a potential $V(\bm{x})=(2|\bm{x}|)^{-1}$ which features a severe point singularity at the origin~$\bm{0}=(0,0)$. Specifically, the energy functional is given by
\begin{align*}
 \E(u)=\frac12\int_\Omega \left( |\nabla u|^2+|\bm{x}|^{-1}|u|^2\right) \dx,
\end{align*}
with $\Omega = \left(\nicefrac{-1}{2},\nicefrac{1}{2}\right)^2$. This experiment has been conducted already in~\cite{MadayMarcati:18}, where the authors obtained an approximated minimal eigenvalue $\lambda_{\mathrm{GS}}=25.934923921168$ for the corresponding GPE~\eqref{eq:GP}; we compare our results to this reference value (divided by~2 on account of \eqref{eq:EVenergy}). We can see from Figure~\ref{fig:alpha1} that the proposed Algorithm~\ref{alg:adaptive} achieves a sequence of energy approximations for the ground state which decays at an (almost) optimal rate.

\subsection{Linear GPE with potential wells} \label{sec:potwells}

\begin{figure}
	\hfill
	\includegraphics[width=0.49\textwidth]{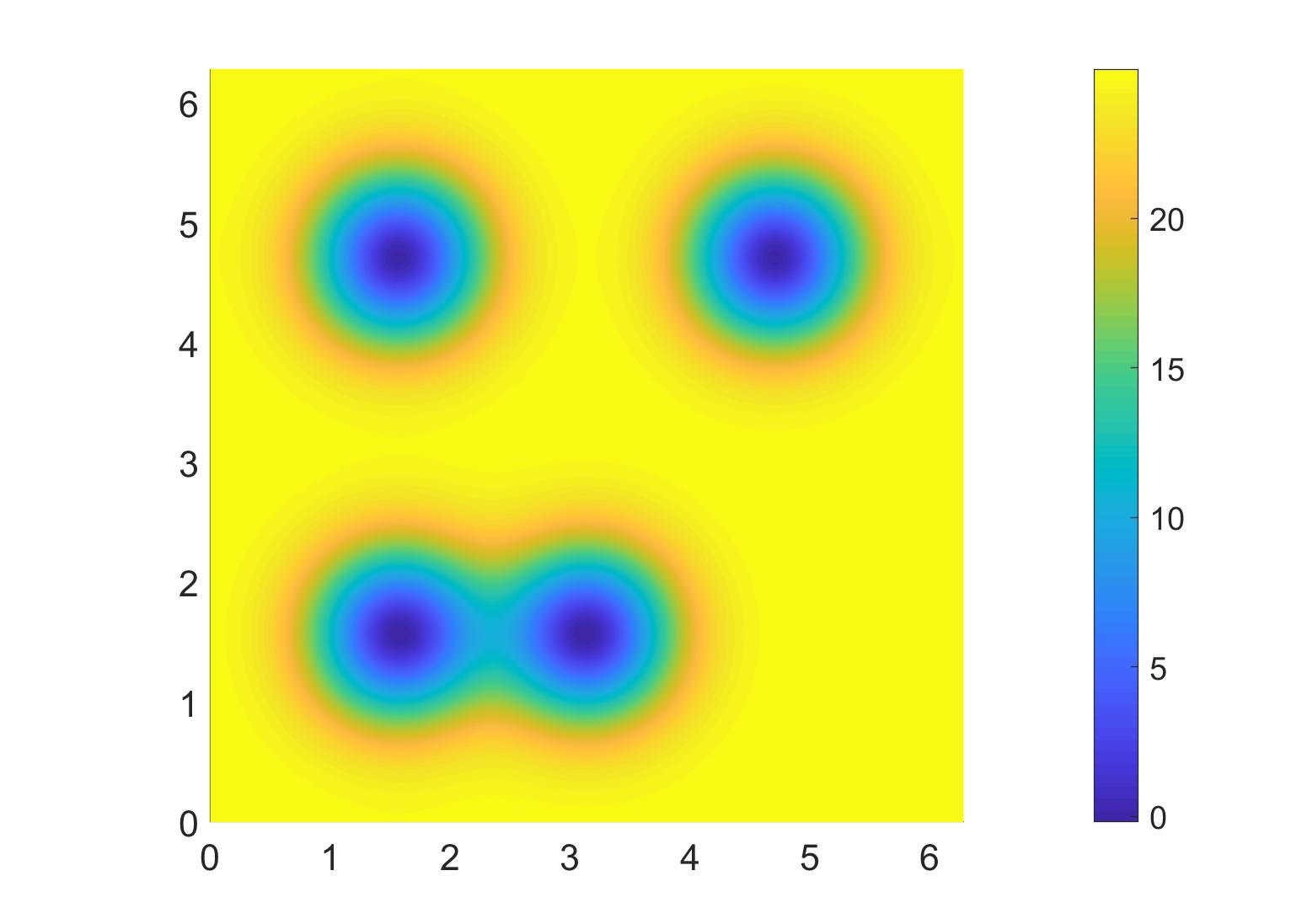}
	\hfill
	\includegraphics[width=0.49\textwidth]{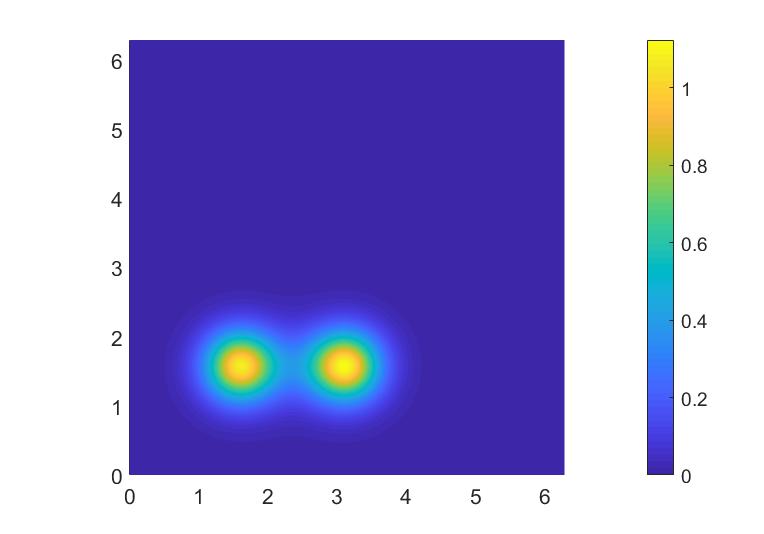}
	\hfill
	\caption{Experiment~\ref{sec:potwells}. Left: Potential function~$V$ consisting of four Gaussians bells. Right: Visualization of the computed ground state.}
	\label{fig:potwells}
\end{figure}

\begin{figure}
 \includegraphics[width=0.49\textwidth]{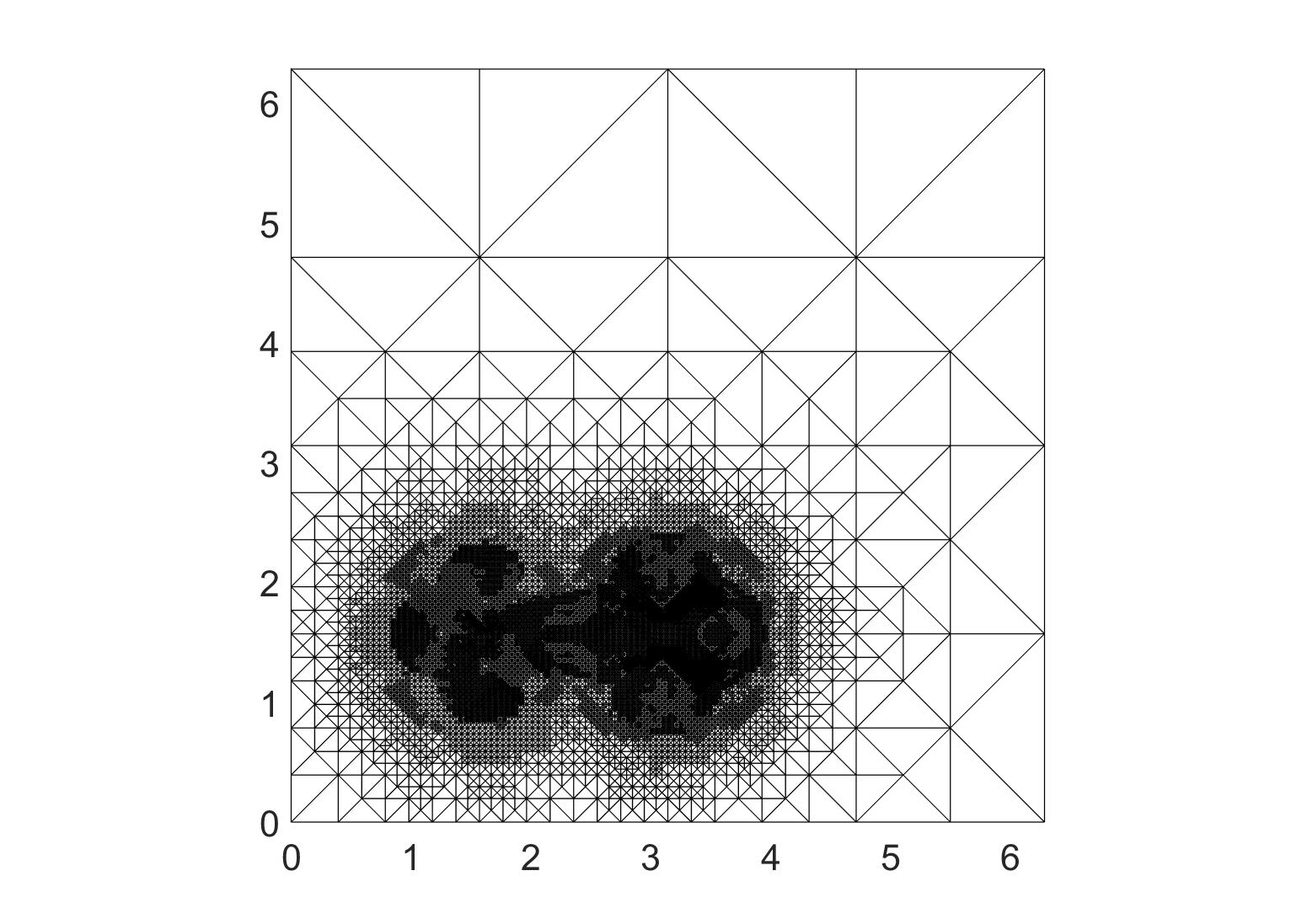}
 \caption{Adaptively refined mesh after 11 refinements.}
 \label{fig:potwells_mesh}
\end{figure}

We test a final example with $\beta=0$. The potential $V$ is given by the sum of four Gaussian bells, see Figure~\ref{fig:potwells} (left). This experiment is borrowed from~\cite[Experiment~4.2]{LinStamm:17}, however, with a (constant) shift such that $V \geq 0$ in the underlying domain~$\Omega=(0,2\pi)^2$. As we can see from Figure~\ref{fig:potwells} (right), the mass of the ground state is mainly concentrated at two adjacent hills in the lower left part of the domain; we note that this is perfectly in line with the results obtained in the paper~\cite{LinStamm:17}. Moreover, the energy-based adaptive mesh refinement has properly resolved the two local hills featured in the ground state $u_{\mathrm{GS}}$, see Figure~\ref{fig:potwells_mesh}.

\begin{figure}
	\hfill
	\includegraphics[width=0.49\textwidth]{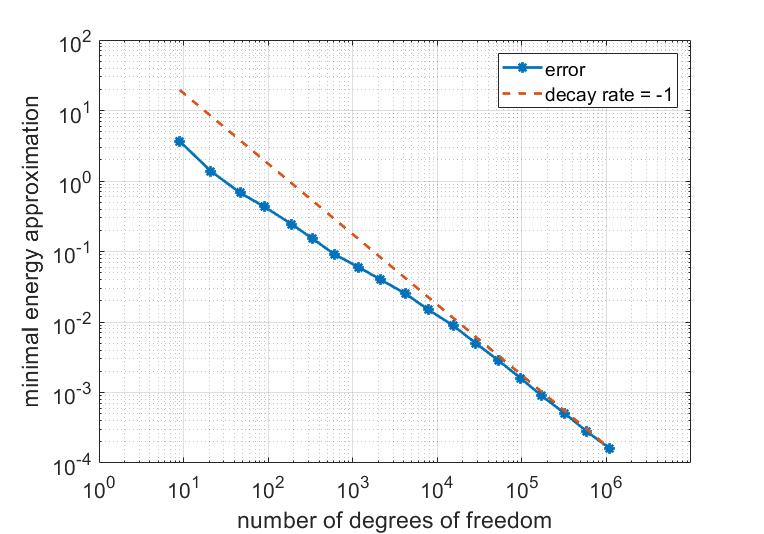}
	\hfill
	\includegraphics[width=0.49\textwidth]{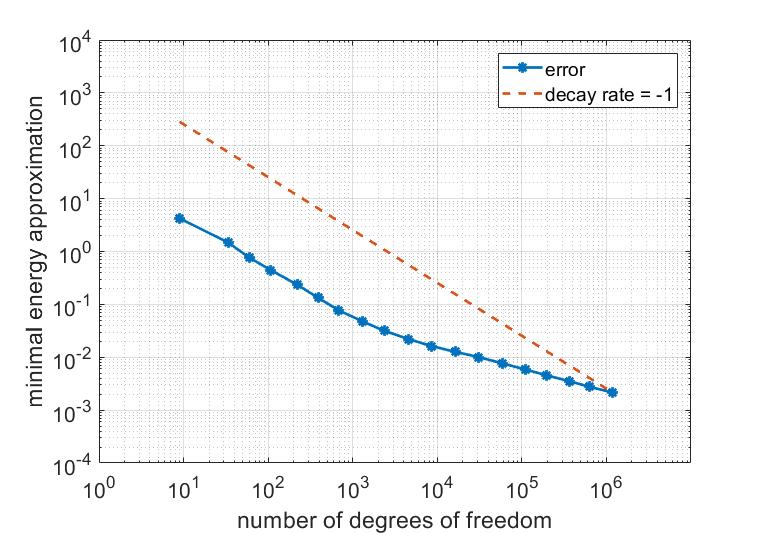}
	\hfill
	\caption{Experiment~\ref{sec:potwells}. Left: Convergence plot for the ground state energy with $\gamma=0.1$. Right: Convergence plot for the ground state energy with $\gamma=0.5$.}
	\label{fig:potwells2}
\end{figure}

\begin{figure}
	\hfill
	\includegraphics[width=0.49\textwidth]{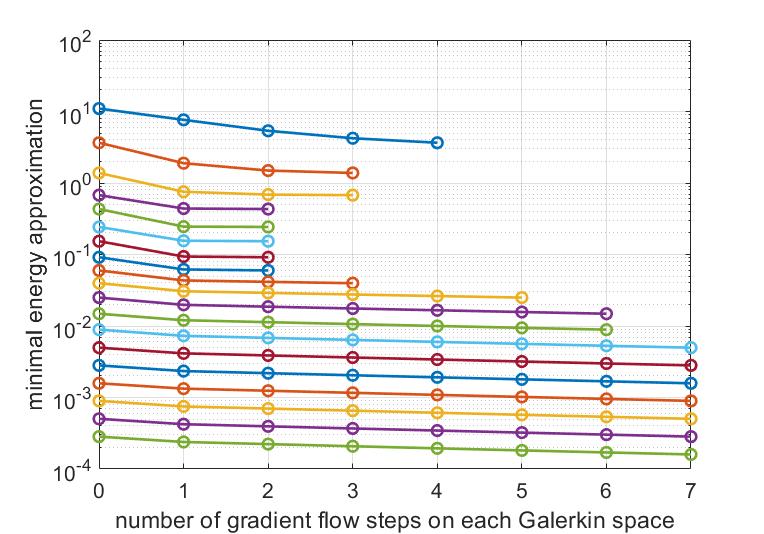}
	\hfill
	\includegraphics[width=0.49\textwidth]{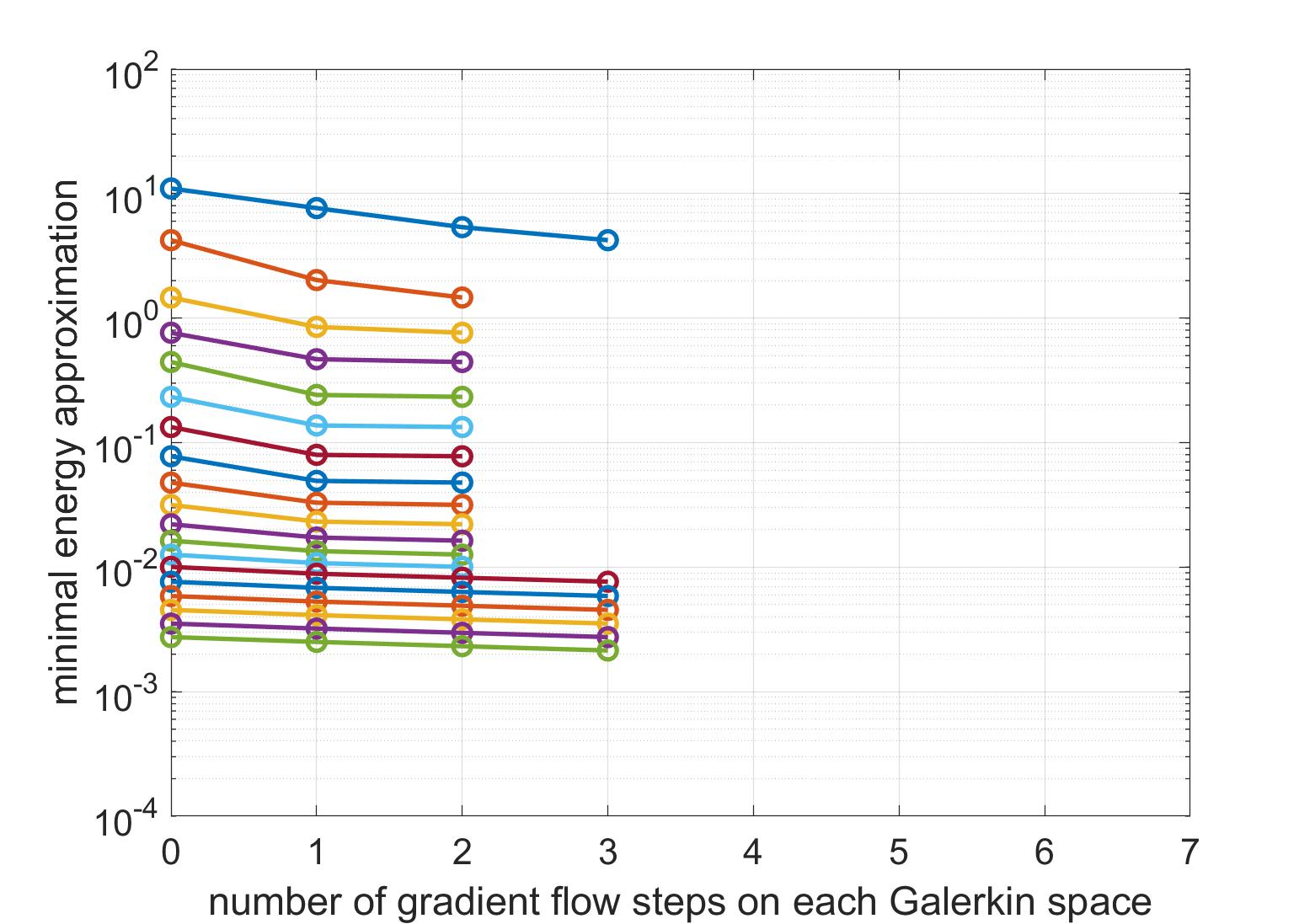}
	\hfill
	\caption{Experiment~\ref{sec:potwells}: Evolution of the minimal energy approximation~$\E(u^n_N)-\E(u_{\mathrm{GS}})$ with respect to the number of gradient flow steps~$n$ and the number of adaptive mesh refinements~$N$, for $\gamma=0.1$ (left) and $\gamma=0.5$ (right): For each Galerkin space $\X_N$ an individual line illustrates the energy decay with respect to the number of gradient flow steps.}
	\label{fig:potwells05_cycle}
\end{figure}

The authors from~\cite{LinStamm:17} have computed an approximated minimal eigenvalue $\lambda_{\mathrm{GS}}=16.6879$, whereby this is the value adapted for our reformulation of the problem. Based on this approximation, we observe an optimal rate of convergence for the energy of the ground state in Figure~\ref{fig:potwells2} (left). We remark, for this example, that the performance of Algorithm~\ref{alg:adaptive}  crucially depends on the choice of the parameter $\gamma$ in Line~6. Indeed, if $\gamma=0.5$ (instead of~$\gamma=0.1$) is selected, then the numerical results exhibit a considerably less favourable asymptotic convergence regime, see Figure~\ref{fig:potwells2} (right). An explanation can be inferred from Figure~\ref{fig:potwells05_cycle}: We see that the choice $\gamma=0.1$ leads to a significantly higher number of gradient flow steps on each Galerkin space $\X_N$, which seems to be essential for the effective numerical solution of this problem.

\subsection{Nonlinear GPE with harmonic confinement potential} \label{sec:beta1000}

\begin{figure}
	\hfill
	\includegraphics[width=0.49\textwidth]{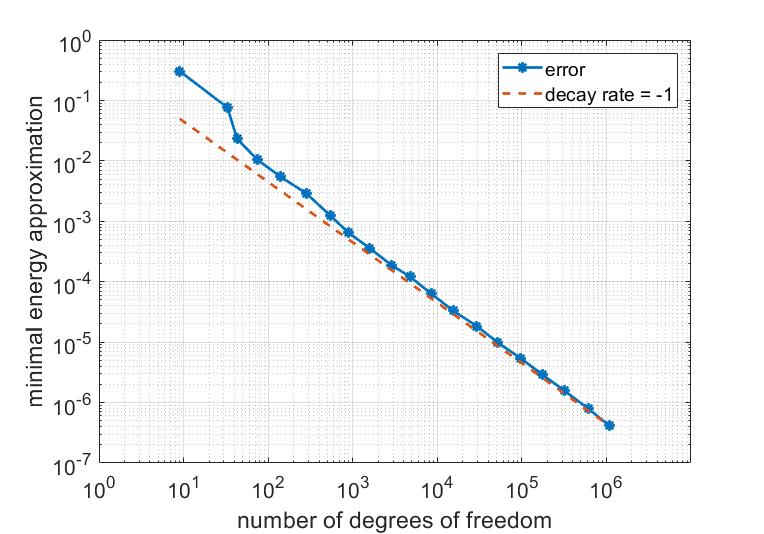}
	\hfill
	\includegraphics[width=0.49\textwidth]{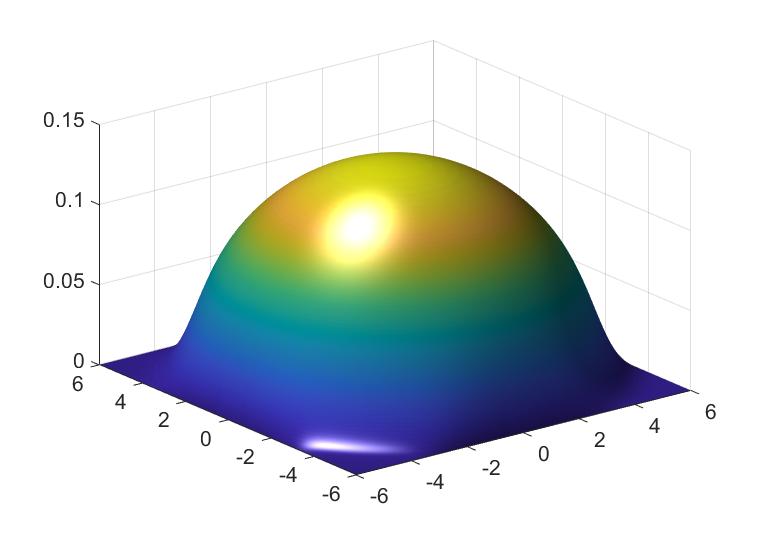}
	\hfill
	\caption{Experiment~\ref{sec:beta1000}. Left: Convergence plot for the ground state energy. Right: Approximated ground state.}
	\label{fig:beta1000}
\end{figure}

We now consider a nonlinear Bose-Einstein condensate, i.e.~$\beta=1000 \gg 0$ in~\eqref{eq:BEC}, and use the smooth potential~$V(\bm{x})=\nicefrac12|\bm{x}|^2$: 
\[
 \E(u)=\frac{1}{2}\int_\Omega\left(  |\nabla u|^2+\left(x^2+y^2\right)|u|^2+1000 |u|^4 \right)\dx; 
\]
the domain is given by $\Omega=(-6,6)^2$. This experiment was conducted previously in~\cite{HenningPeterseim:18}, where an approximation $\E(u_\mathrm{GS}) \approx 11.9860647$ for the energy of the ground state has been documented. Based on the adaptive Algorithm~\ref{alg:adaptive} presented in this work, a smaller value for the  ground state energy has been computed; we suppose that this (improved) approximation results from the adaptive (and thereby more effective) refinement of the meshes. To be specific, we have obtained the approximation $\E(u_\mathrm{GS}) \approx \underline{11.98605}121\dotsc$ based on an adaptively refined mesh with $\mathcal{O}(10^7)$ degrees of freedom; here, the underlined digits are stable (i.e.~the computations indicate that they do not change any more as the iterations continue). Using this as reference value, we obtain optimal convergence for the approximation of the energy of the ground state, see Figure~\ref{fig:beta1000}.

In order to study the dependence of the computational work on the strength of the nonlinearity, we have performed the same experiment for different values of $\beta$. 
As can be seen from Table~\ref{table:beta1}, our algorithm does not seem to be sensitive with respect this parameter. Indeed, after an initial phase, the number of GF iterations on a given mesh is independent of~$\beta$. Moreover, for $\gamma=0.5$, this holds true already from the start, i.e. even for the initial mesh (not displayed here). We remark that analogous results have been observed recently in~\cite{ZhangXuXie:19} for a multigrid method.

\begin{center}
\small{
  \begin{tabular}{c|ccccccccccccccc} 
    \toprule
    $\beta$ & $\X_0$ &$\X_1$ & $\X_2$ & $\X_3$ &$\X_4$ &$\X_5$&$\X_6$&$\X_7$&$\X_8$&$\X_9$&$\X_{10}$&$\X_{11}$&$\X_{12}$&$\X_{13}$&$\X_{14}$ \\ \midrule
    0 & 2 & 2 & 2 &2 &2&2&2&2&2&2&2&2&2&2&2 \\ 
    200 & 2 & 3 & 3 & 3 &2&2&2&2&2&2&2&2&2&2&2 \\ 
    400 & 3 & 3 & 3 & 3 &3&3&3&2&2&2&2&2&2&2&2 \\ 
    600 & 3 & 3 & 3 & 4 &3&3&2&2&2&2&2&2&2&2&2 \\ 
    800 & 3 & 3 & 4 & 4 &3&3&3&2&2&2&2&2&2&2&2 \\ 
	\bottomrule
    \end{tabular}
    }
      \captionof{table}{Experiment~\ref{sec:beta1000}. Number of GF iterations on the first 15 finite element spaces for different values of $\beta$.}
          \label{table:beta1}
\end{center}

\subsection{Nonlinear GPE with optical lattice potential} \label{sec:complicatedV}

\begin{figure}
	\hfill
	\includegraphics[width=0.49\textwidth]{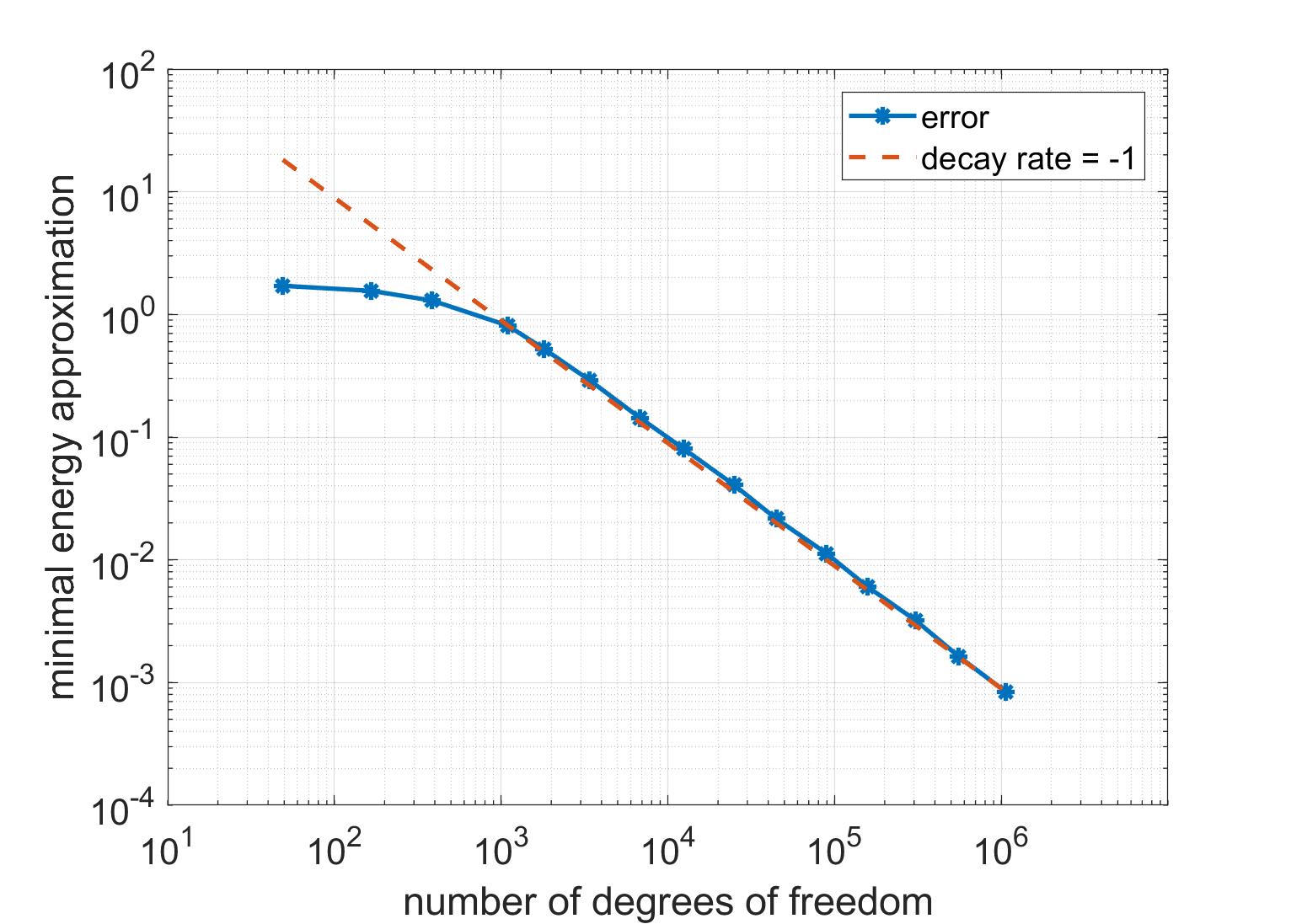}
	\hfill
	\includegraphics[width=0.49\textwidth]{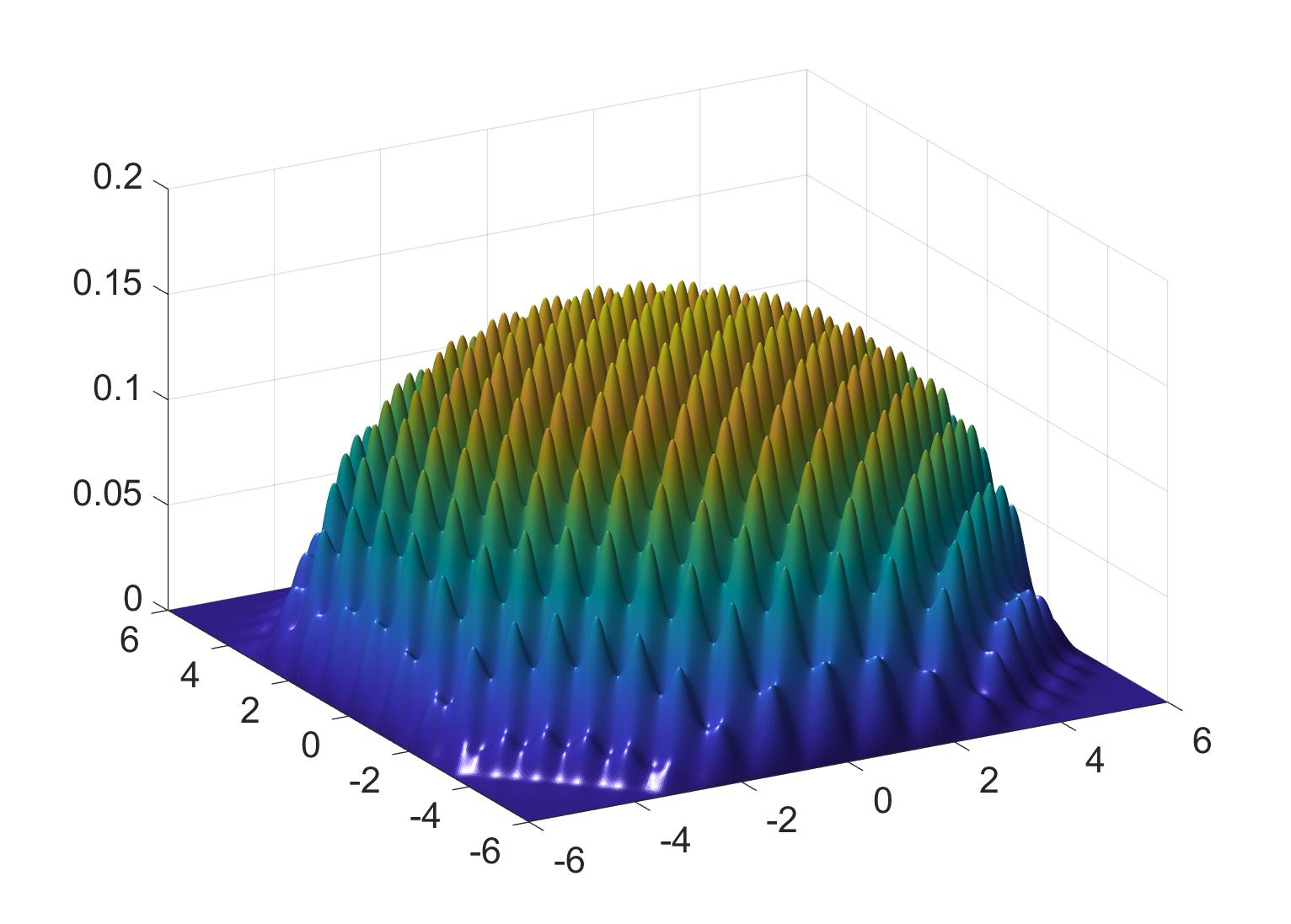}
	\hfill
	\caption{Experiment~\ref{sec:complicatedV}. Left: Convergence plot for the ground state energy. Right: Approximated ground state.}
	\label{fig:complicatedV}
\end{figure}

\begin{figure}
	\hfill
	\includegraphics[width=0.49\textwidth]{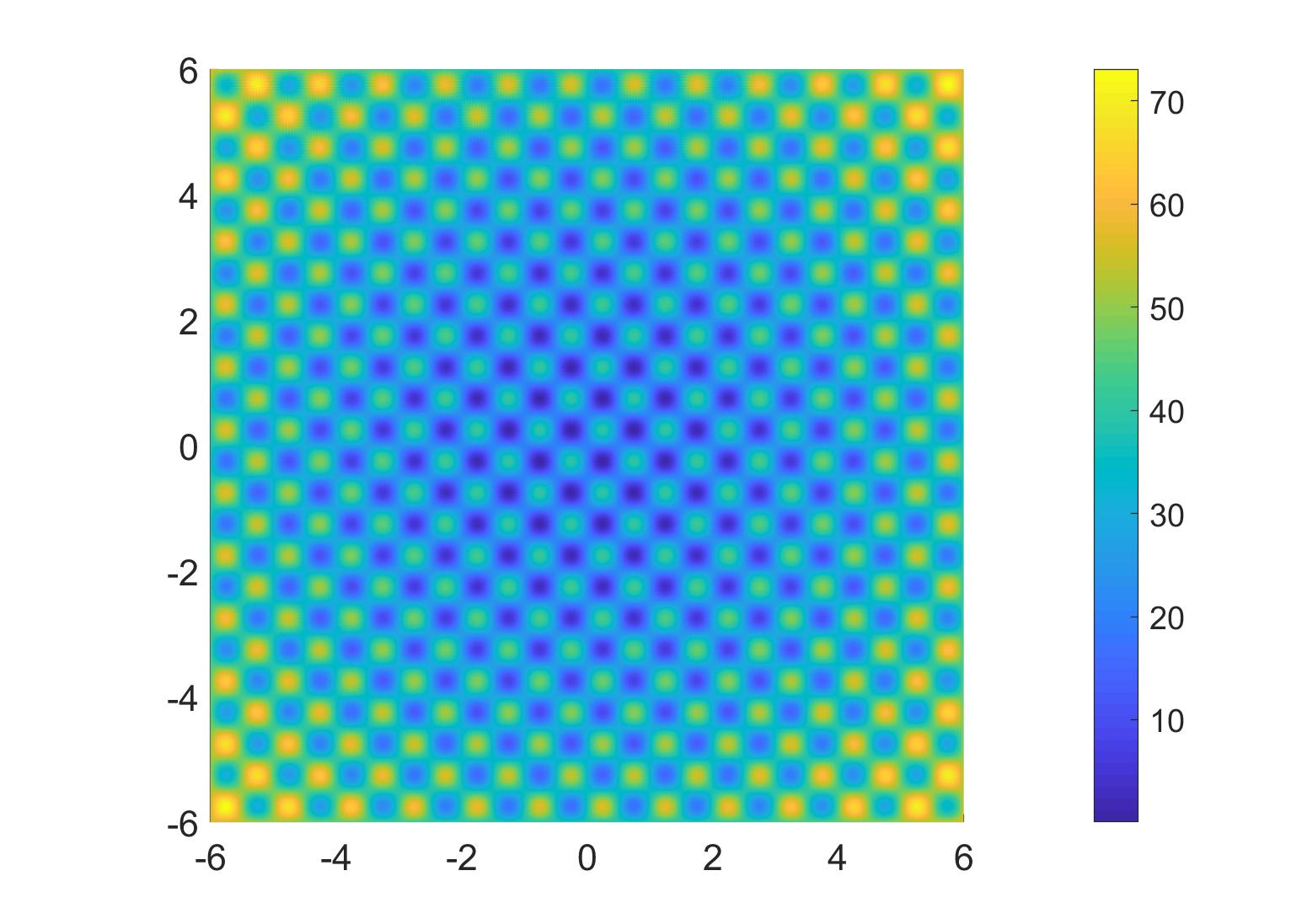}
	\hfill
	\includegraphics[width=0.49\textwidth]{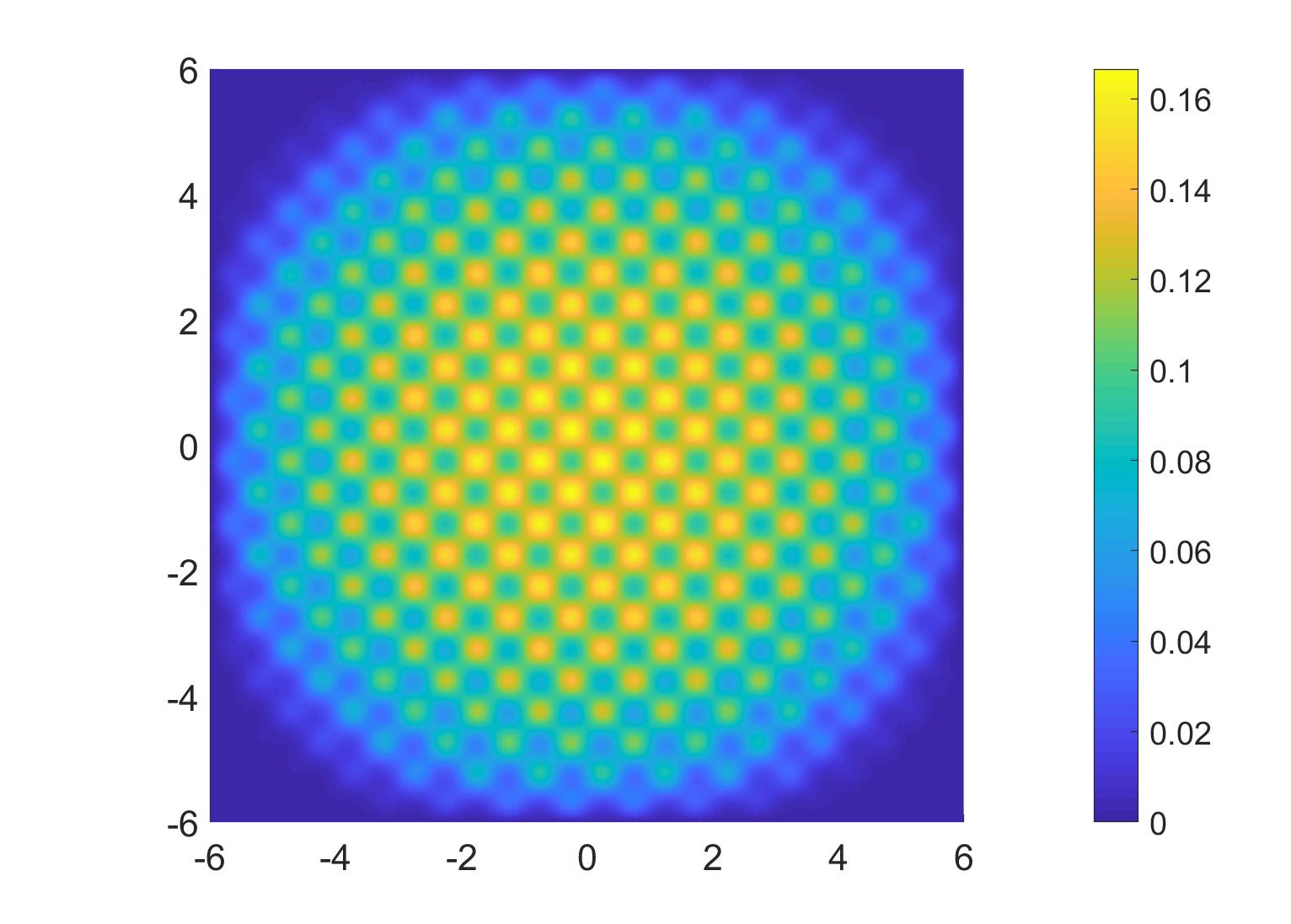}
	\hfill
	\caption{Experiment~\ref{sec:complicatedV}. Left: Contour plot of the potential $V$. Right: Contour plot of the ground state.}
	\label{fig:complicatedV2}
\end{figure}

As before, we choose $\Omega=(-6,6)^2$ and $\beta=1000$, with an oscillating potential function~$V$, see Figure~\ref{fig:complicatedV2} (left) for its contour plot. More precisely, the energy functional is given by
\[
  \E(u)=\frac{1}{2}\int_\Omega \left( |\nabla u|^2+2\left(\frac{|\bm{x}|^2}{2}+20+20 \sin(2\pi x) \sin(2\pi y)\right)|u|^2+1000 |u|^4\right) \dx.
\]
This experiment was also considered in~\cite{HenningPeterseim:18} with an asserted approximation $\E(u_\mathrm{GS}) \approx 30.40965$ of the ground state energy. Based on a sufficiently fine uniform initial mesh (with $128$ elements), our algorithm yields the approximation $\E(u_\mathrm{GS}) \approx \underline{30.387}533\dotsc$ for $\mathcal{O}(10^7)$ degrees of freedom. In Figure~\ref{fig:complicatedV} (left) we have depicted the error for the approximations of the ground state energy with respect to our reference value. This plot indicates an asymptotically optimal rate of convergence of Algorithm~\ref{alg:adaptive} for the given problem. 

\begin{figure}
	\hfill
	\includegraphics[width=0.49\textwidth]{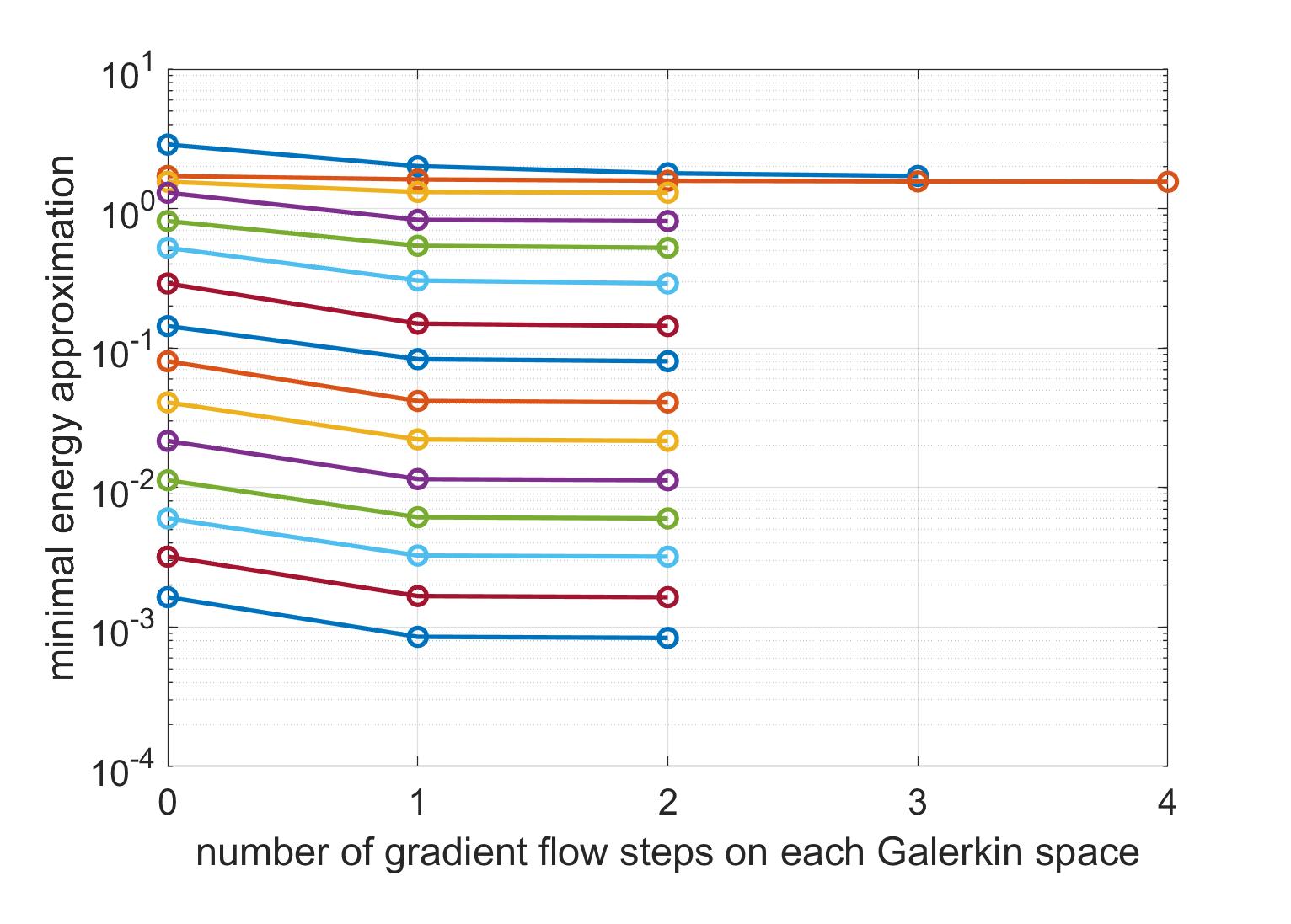}
	\hfill
	\includegraphics[width=0.49\textwidth]{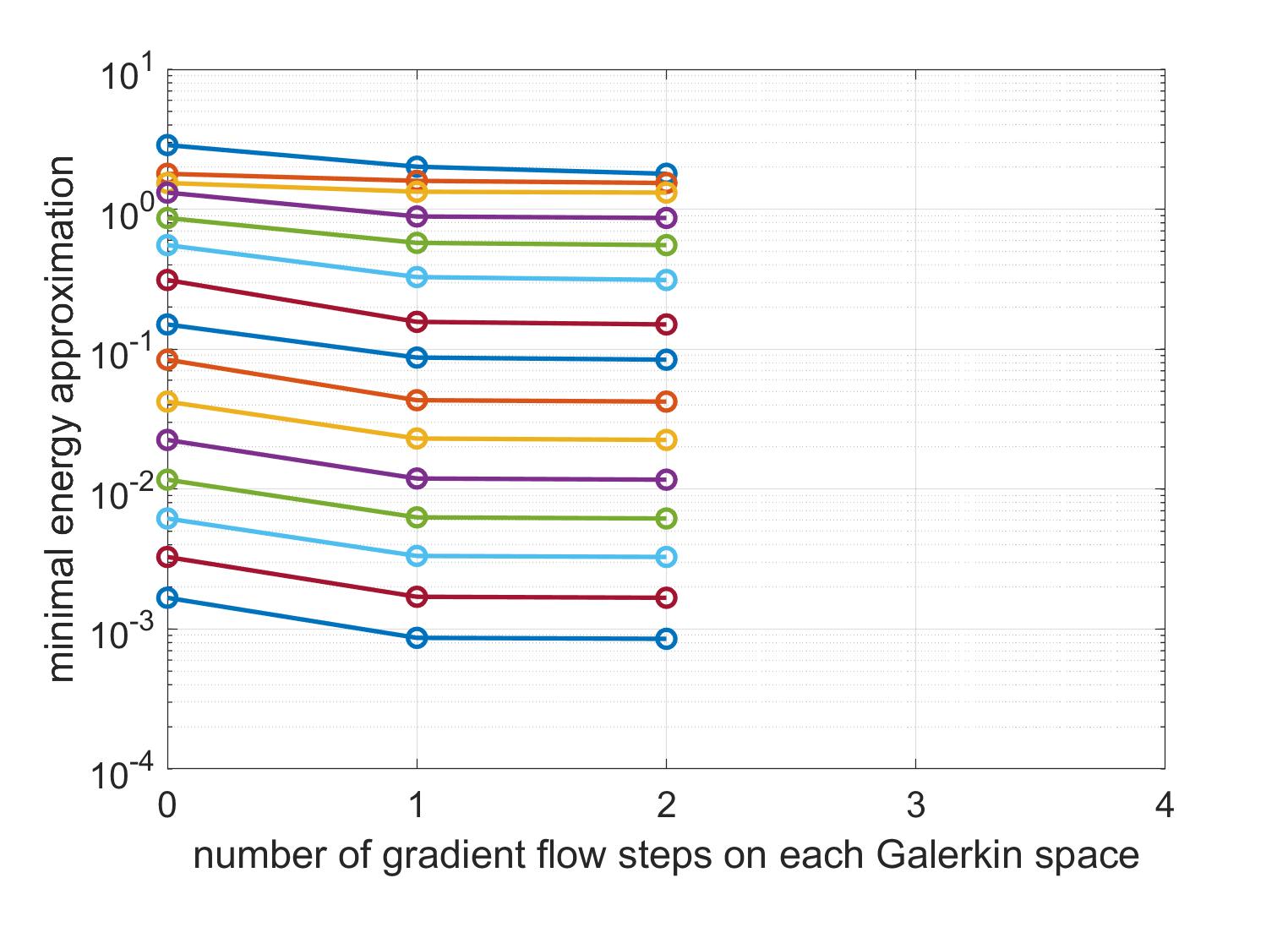}
	\hfill
	\caption{Experiment~\ref{sec:complicatedV}. Evolution of the minimal energy approximation~$\E(u^n_N)-\E(u_{\mathrm{GS}})$ with respect to the number of gradient flow steps~$n$ and the number of adaptive mesh refinements~$N$, for $\gamma=0.1$ (left) and $\gamma=0.5$ (right).}
	\label{fig:NonlinearComplicated_cycle}
\end{figure}

As in Example~\ref{sec:potwells}, we have also run this experiment for both values $\gamma=0.1$ and~$\gamma=0.5$.  In contrast to the previous test, we observe no considerable difference in the performance of the corresponding computations. Indeed, except for the first two meshes, the number of gradient flow steps on each Galerkin space is the same for both cases, see Figure~\ref{fig:NonlinearComplicated_cycle}.

\subsection{Nonlinear energy functional with a nonsymmetric potential $V$} \label{sec:nonsymmetric}

\begin{figure}
	\hfill
	\includegraphics[width=0.49\textwidth]{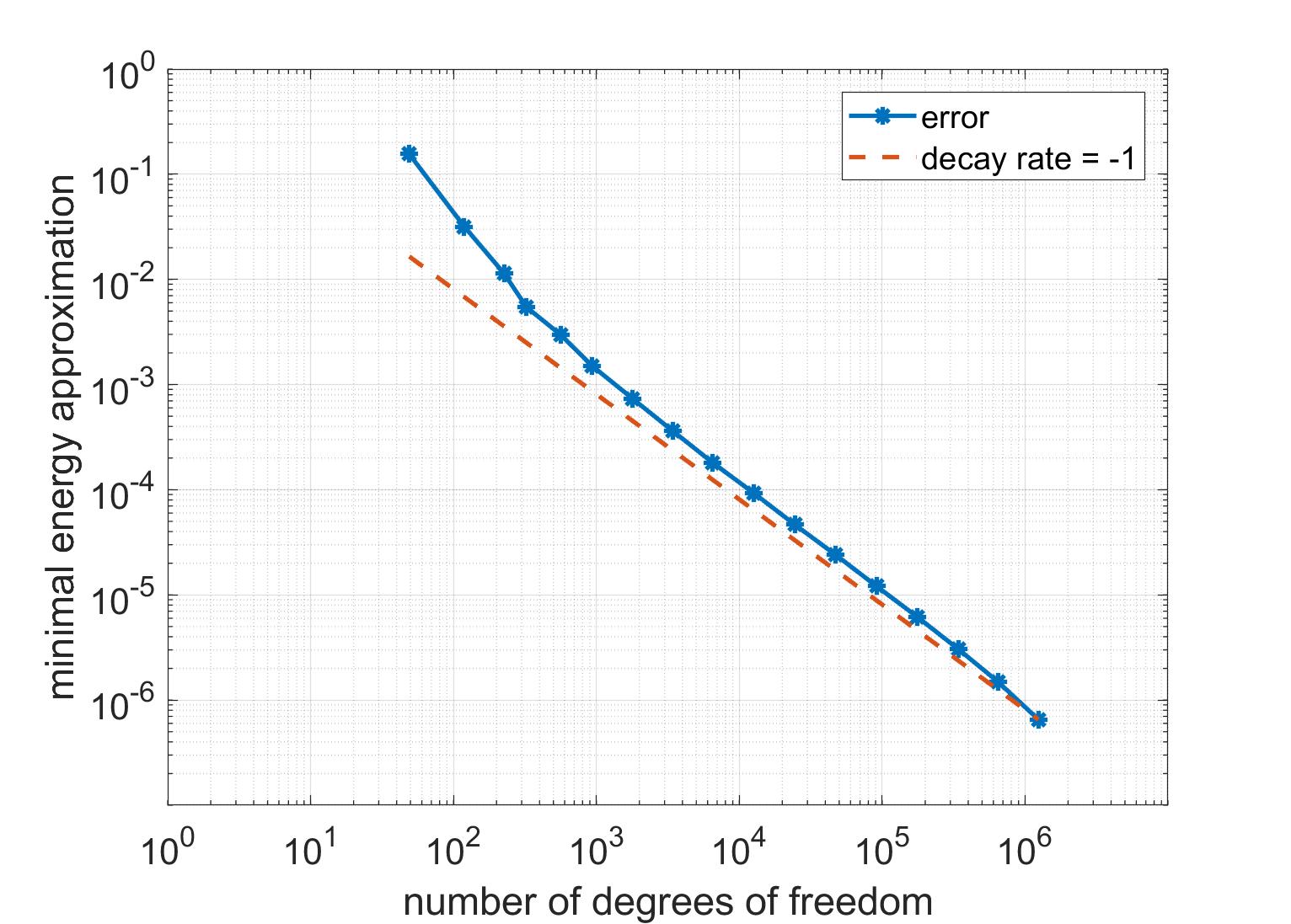}
	\hfill
	\includegraphics[width=0.49\textwidth]{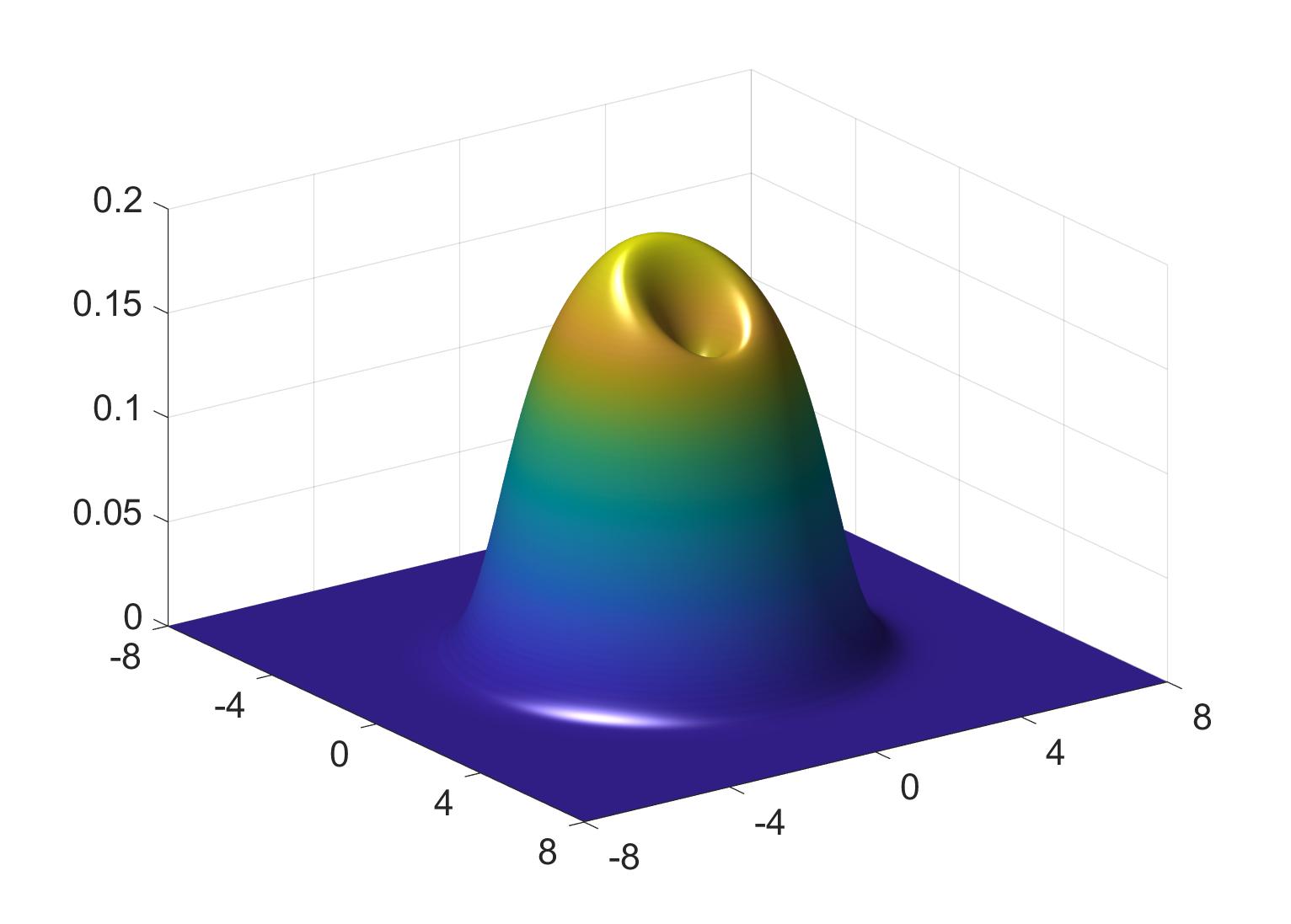}
	\hfill
	\caption{Experiment~\ref{sec:nonsymmetric}. Left: Convergence plot for the ground state energy. Right: Approximated ground state.}
	\label{fig:nonsymmetric_convergence}
\end{figure}

Finally, we revisit the test problem~\cite[Example~4.3.II]{BaoDu:04}:\
\[
  \E(u)= \frac{1}{2} \int_\Omega \left( |\nabla u|^2 +\left(|\bm{x}|^2+8 \exp\left(-((x-1)^2+y^2)\right)\right)|u|^2+200 |u|^4\right) \dx,
\]
with the symmetric domain $\Omega=(-8,8)^2$. In~\cite{BaoDu:04} an approximated value of $\E(u_{\mathrm{GS}}) \approx 5.8507$ for the energy of the ground state $u_{\mathrm{GS}}$ has been obtained. Like in the previous experiments, for an adaptively refined mesh with $\mathcal{O}(10^7)$ degrees of freedom, a smaller value of $\E(u_{\mathrm{GS}}) \approx \underline{5.85058}738\dotsc$ is obtained. Applying this value as reference ground state energy, we observe optimal convergence in Figure~\ref{fig:nonsymmetric_convergence}.
Furthermore, following the previous Experiment~\ref{sec:beta1000}, we have run this test problem for different values of $\beta$, with similar results, see Table~\ref{table:beta2}.

\begin{center}
\small{
\begin{tabular}{ c|ccccccccccccccc }
    \toprule
    $\beta$ & $\X_0$ &$\X_1$ & $\X_2$ & $\X_3$ &$\X_4$ &$\X_5$&$\X_6$&$\X_7$&$\X_8$&$\X_9$&$\X_{10}$&$\X_{11}$&$\X_{12}$&$\X_{13}$&$\X_{14}$ \\ \midrule
    0 & 3  &   3  &   2  &   2  &   2  &   2  &   2  &   2  &   2  &   2  &   2  &   2 & 2 & 2 & 2 \\ 
    200 & 3  &   3  &   3  &   3  &   2  &   2  &   2  &   2  &   2  &   2  &   2  &   2 & 2 & 2 & 2 \\ 
    400 & 3  &   3  &   3  &   3  &   2  &   2  &   2  &   2  &   2  &   2  &   2  &   2 & 2 & 2 & 2 \\ 
    600 & 3  &   3  &   4  &   3  &   2  &   2  &   2  &   2  &   2  &   2  &   2  &   2 & 2 & 2 & 2\\ 
    800 & 3  &   4  &   4  &   3  &   3  &   2  &   2  &   2  &   2  &   2  &   2  &   2 & 2 & 2 & 2 \\ \bottomrule

    \end{tabular}
    }
      \captionof{table}{Experiment~\ref{sec:nonsymmetric}. Number of GF iterations on the first 15 finite element spaces for different values of $\beta$.}
           \label{table:beta2}
\end{center}

\section{Conclusions}\label{sec:conclusions}

In this work, we have considered a computational procedure for the numerical approximation of the ground state and its associated energy of the Gross-Pitaevskii equation, which applies an effective interplay of a gradient flow iteration method and adaptive mesh refinements. Both of these techniques rely on energy minimization and guaranteed energy reduction. Thereby, they are based on the underlying structure of the problem at hand in a very natural way. Our scheme is fairly simple to implement and, for the test problems presented here, exhibits either optimal or close to optimal convergence rates for the approximation of the ground state energy. Moreover, in our experiments, the effect of the parameter $\beta$ (steering the strength of the nonlinearity)  on the computational work seems to be negligible.

\bibliographystyle{amsplain}
\bibliography{references,GPrefs}

\end{document}